\documentclass[11pt]{article}

\usepackage{fullpage}
\usepackage{graphicx}
\usepackage{fancyvrb}
\usepackage{xcolor}
\usepackage{amsmath}
\usepackage{enumerate}
\usepackage{amssymb}
\usepackage{amsfonts}
\usepackage{natbib} 

\usepackage{float}
\usepackage{bigints}

\def\rr{{\mathbb R}}
\def\ee{{\mathbb E}}
\def\bX{{\mathbf X}}
\def\bx{{\mathbf x}}
\def\by{{\mathbf y}}
\def\bY{{\mathbf Y}}
\def\b0{{\mathbf 0}}

\newcommand{\Linf}{$L^\infty$}
\newcommand{\Ltwo}{$L^2$}

\usepackage{ulem}

\begin{document} 

\title{\bf Numerical studies of space filling designs:
optimization of Latin Hypercube Samples
and subprojection properties} 

\author{\bf G. Damblin, M. Couplet and B. Iooss \\
EDF R\&D, 6 Quai Watier, F-78401, Chatou, France 
}
\date{}
\maketitle 

\begin{center}
  Submitted to: {\it Journal of Simulation}\\
for the special issue ``Input \& Output Analysis for Simulation''

  \vspace{0.2cm}
  Correspondance: B. Iooss ; Email: bertrand.iooss@edf.fr\\
  Phone: +33-1-30877969 ; Fax: +33-1-30878213
\end{center}

\renewcommand{\baselinestretch}{1.5} 


\abstract{
Quantitative assessment of the uncertainties tainting the results of computer simulations is nowadays a major topic of interest in both industrial and scientific communities.
One of the key issues in such studies is to get information about the output when the numerical simulations are expensive to run. 
This paper considers the problem of exploring the whole space of variations of the computer model input variables in the context of a large dimensional exploration space.
Various properties of space filling designs are justified: interpoint-distance, discrepancy, minimum spanning tree criteria.
A specific class of design, the optimized Latin Hypercube Sample, is considered.
Several optimization algorithms, coming from the literature, are studied in terms of convergence speed, robustness to subprojection and space filling properties of the resulting design.
Some recommendations for building such designs are given.
Finally, another contribution of this paper is the deep analysis of the space
filling properties of the design 2D-subprojections.
}



\section{Introduction}
\label{sec:intro}

Many computer codes, for instance simulating physical
phenomena and industrial systems, are too time expensive to be directly
used to perform uncertainty, sensitivity, optimization or
robustness analyses \citep{derdev08}. 
A widely accepted method to circumvent this problem consists in replacing
such computer models by cpu time
inexpensive mathematical functions, called metamodels \citep{klesar00},
built from a limited number of
simulation outputs.
Some commonly used
metamodels are: polynomials, splines, generalized linear models, or learning
statistical models like neural networks, regression trees, support vector
machines and Gaussian process models \citep{fanli06}. 
In particular,
the efficiency of Gaussian process
modelling has been proved for instance by \cite{sacwel89,sanwil03,marioo07}.
It extends the kriging principles of geostatistics to computer experiments by
considering that the code responses are correlated
according to the relative locations of the corresponding input variables.

A necessary condition to successful metamodelling is to explore the whole space
$\mathcal{X} \subset \rr^{d}$ of the
 input variables $\bX \in \mathcal{X}$
(called the inputs) 
in order to  capture the non linear
behaviour of some output variables
$\bY=G(\bX) \in \rr^q$ (where $G$ refers to the computer
code). This step, often called the Design of Computer Experiments (DoCE),
is the subject of this paper.
In many industrial applications, we are faced with the harsh
problem of the high dimensionality of the space
$\mathcal{X}$ to explore (several tens of inputs).
For this purpose, Sample Random Sampling (SRS), that is standard
Monte-Carlo, can be used. An important advantage of SRS is that the linear
variability of estimates of quantities such as mean output behaviour
is independent of $d$, the dimension of the input space.

However, some authors 
\citep{simlin01,fanli06} have shown that
the Space Filling Designs (SFD) are
well suited to this initial exploration objective.
A SFD aims at obtaining the best
coverage of the space of the inputs,
while SRS does not ensure this task.
Moreover, the SFD approach appears natural in first exploratory phase,
when very few pieces of information are available about the numerical
model, or if the design is expected to serve different objectives (for example,
providing a metamodel usable for several uncertainty quantifications
relying on different hypothesis about the uncertainty on the inputs).
However, the class of SFD is diverse including the
well-known Latin Hypercubes Samples (LHS) \footnote{
In the following, LHS may refer to Latin Hypercube Sampling as well.},
Orthogonal Arrays (OA), point process designs \citep{frabay08}, minimax and
maximin designs each based on a distance criterion \citep{johmoo90,koeowe96},
maximum entropy designs \citep{schwyn87}. Another kind of SFD is the one
of Uniform Designs (UD) whose many interests for computer experiments have been already
discussed \citep{fan03}. They are based upon some measures of
uniformity called discrepancies \citep{nider87}. Here, our purpose is to shed
new light on the practical issue of building a SFD.
 
In the following, $\mathcal{X}$ is assumed to be ${[0;1]}^d$, up to
a bijection. Such a bijection is never unique and different ones generally
lead to non equivalent ways of filling the space of the inputs.
In practice, during an exploratory phase, only pragmatic answers can be
given to the question of the choice of the bijection: maximum and
minimum bounds are generally given to each scalar input so that $\mathcal{X}$
is an hypercube which can be mapped over ${[0;1]}^d$ through a linear
transformation.
It can be noticed that, if there is a sufficient reason to do so,
it is always possible
to apply simple changes of input variables if it seems relevant
(\textit{e.g.} considering the input $z = \exp(x) \in [\exp(a),\exp(b)]$
instead of $x \in [a, b]$).

In what follows, we focus on the
``homogeneous'' filling of ${[0;1]}^d$ and
the main question addressed is how to build or to select a DoCE 
of a given (and small) size $N$
($N \sim 100$, typically)
within ${[0;1]}^d$ (where $d > 10$, typically).
We keep also in mind the well-known and empirical
relation $N \sim 10 d$
\citep{loesac09,marioo08} which gives the approximative
minimum number of computations needed to get an accurate metamodel.

A first consideration is to obtain the best global coverage rate of ${[0;1]}^d$.
In the following, quality of filling is improved by optimizing discrepancies
or point-distance criteria and assessed by robust geometric criteria based on
the Minimum Spanning Tree (MST) \citep{dusras86}.
When building SFD upon a particular criterion\footnote{
Which includes boolean ones, that is properties like
''having a Latin Hypercube structure``.},
a natural prerequisite is coordinate rotation invariance, that is
invariance of the value of the criterion under exchanges of
coordinates. Criteria applied hereafter are coordinate rotation invariant.
An important prescription is to properly cover the variation domain
of each scalar input. Indeed, it often happens that,
among a large number of inputs, only a small
one is active, that is significantly impacts the outputs
(sparsity principle).
Then, in order to avoid useless
computations (different values for inactive inputs but same values for active
ones), we have to ensure that all the values for each
input are different, which can be achieve
by using LHS.
A last important property of a SFD is its robustness to projection over subspaces.
This property is particularly studied in this paper and
the corresponding results can be
regarded as the main contributions.
Literature about the application of the physical experimental design theory 
shows that, in most of the practical cases, effects
of small degree (involving few factors, that is few inputs)
dominates effects of greater degree.
Therefore, it seems judicious to favour a SFD whose subprojections 
offer some good coverages of the low-dimensional subspaces. 
A first view, which is adopted here,
is to explore two-dimensional (2D) subprojections.

Section 2 describes two industrial examples in order to
motivate our concerns about SFD.
Section 3 gives a review about
coverage criteria and natures of SFD
studied in the next section.
The preceding considerations lead us to focus our attention on
optimized LHS for which
various optimization algorithms have been previously
proposed (main works are \cite{mormit95} and \cite{jinche05}).
We adopt a numerical approach to compare the performance of different LHS, in
terms of their interpoint-distance and \Ltwo-discrepancies.
Section 4 focuses on their 2D subprojection properties and
numerical tests support some recommendations.
Section 5 synthesizes the work.

\section{Motivating examples}
\label{sec:examples}

\subsection{Nuclear safety simulation studies}\label{sec:nuclear}

Assessing the performance of nuclear power
plants during accidental transient conditions has been the main purpose of
thermal-hydraulic safety research for decades.
Sophisticated computer codes have been developed and
are now widely used. They can calculate time trends
of any variable of interest during any transient in Light Water Reactors (LWR).
However, the reliability of the predictions cannot be
evaluated directly due to the lack of suitable
measurements in plants. The capabilities of the codes can consequently only be
assessed by comparison of calculations with
experimental data recorded in small-scale facilities.
Due to this
difficulty, but also the ``best-estimate'' feature of the
codes quoted above, 
uncertainty quantification should be performed when using them.
In addition to uncertainty quantification,
sensitivity analysis is often carried
out in order to identify the main contributors to uncertainty.

\begin{figure}[!ht]
$$\includegraphics[scale=0.65]{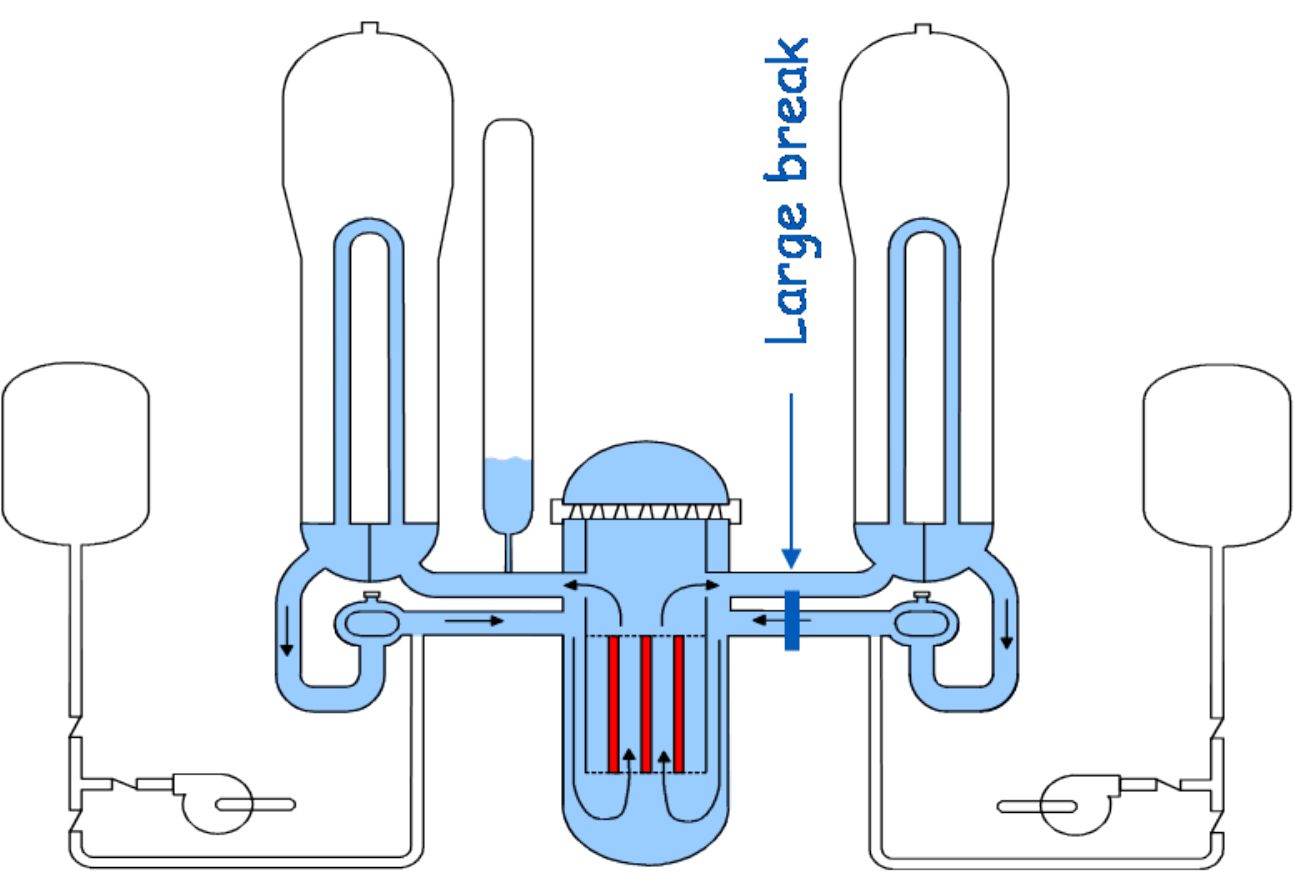}$$
\caption{Illustration of a large-break loss of primary coolant accident
on a nuclear Pressurized Water Reactor (a particular but common
type of LWR).}\label{bemuse}
\end{figure}

\begin{figure}[!ht]
$$\includegraphics[scale=0.55]{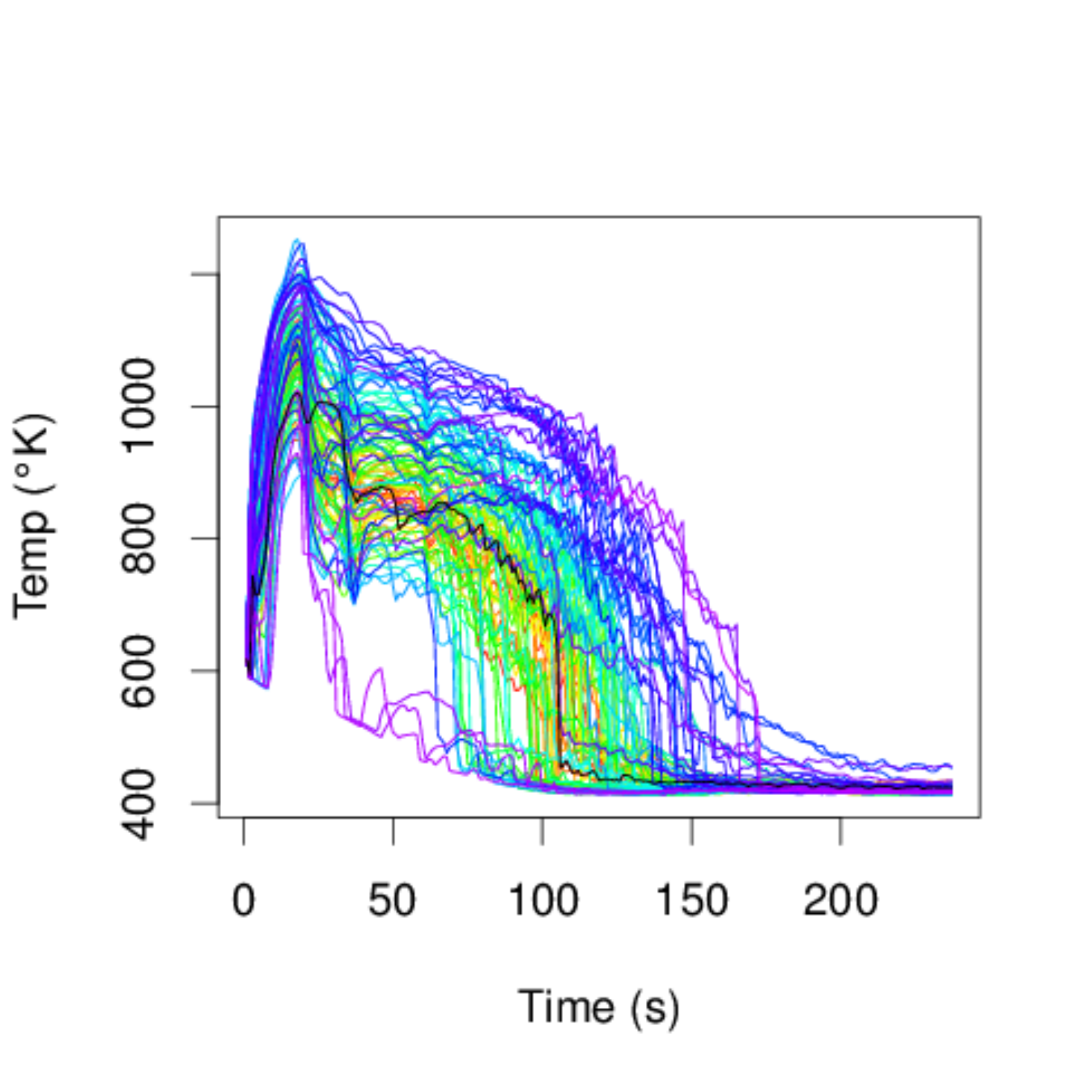}$$

\vspace{-1cm}
\caption{$100$ output curves of the cladding
temperature in function of time from CATHARE2.}
\label{sorties_bemuse}
\end{figure}

Those thermal-hydraulic codes enable, for example, to simulate
a large-break loss of primary coolant accident (see Figure \ref{bemuse}).
This scenario is part of the Benchmark for Uncertainty Analysis in
Best-Estimate Modelling  for Design, Operation and Safety Analysis of Light
Water Reactors \citep{decbaz08} proposed by the Nuclear Energy Agency of the
Organisation for Economic Co-operation and Development (OCDE/NEA).
It has been implemented on the French computer code CATHARE2 developed at the
Commissariat \`a l'Energie Atomique (CEA).
Figure \ref{sorties_bemuse} illustrates $100$ Monte Carlo simulations (by
randomly varying the inputs of the accidental
scenario), given by CATHARE2, of the cladding temperature in function of
time, whose first peak is the main output of interest in safety studies.

Severe difficulties arise when carrying out a sensitivity analysis or an uncertainty quantification
involving CATHARE2:
\begin{itemize}
\item Physical models involve complex phenomena (non linear and
subject to threshold effects), with strong interactions between
inputs. A first objective is to detect these interactions.
Another one is to fully explore combinations of the input to obtain a good idea of the possible
transient curves \citep{auddec12}.
\item Computer codes are cpu time expensive: no more than several hundreds of
simulations can be performed.
\item Numerical models take as inputs a large number of uncertain variables
($d=50$, typically): physical laws essentially, but also initial conditions,
material properties and geometrical parameters.
Truncated normal or log-normal
distributions are given to them. Such a number of inputs is
extremely large for the metamodelling problem. 
\item The first peak of the cladding
temperature can be related to rare events: problems turn to the estimation of a quantile
\citep{cangar08} or the probability that the output exceeds a threshold
\citep{marmar12}.
\end{itemize}

All of these four difficulties underline the fact that 
great care is required to define an effective DoCE over the CATHARE2 input space.
The high dimension of the input space remains  a challenge
for building a SFD with good subprojection properties.

\subsection{Prey-predator simulation chain}\label{sec:prey}

In ecological effect assessments, risks imputable to chemicals are usually
estimated by extrapolation of single-species toxicity tests.
With increasing awareness of the importance of indirect effects and keeping in
mind limitations of experimental tools, a number of ecological food-web models
have been developed.
Such models are based on a large set of bioenergetic equations describing the
fundamental growth of each population, taking into account
grazing and predator-prey interactions, as well as influence of abiotic
factors like temperature, light and nutrients.
They can be used for
several purposes, for instance:
\begin{itemize}
\item to test various contamination scenarios or recovery capacity of
contaminated ecosystem, 
\item to quantify the important sources of uncertainty and knowledge
gaps for which additional data are needed, and to identify the most
influential parameters of population-level impacts,
\item to optimize the design of field or mesocosm tests by identifying
the appropriate type, scale, frequency and duration of monitoring.
\end{itemize}

\begin{figure}[!ht]
$$\includegraphics[scale=0.8]{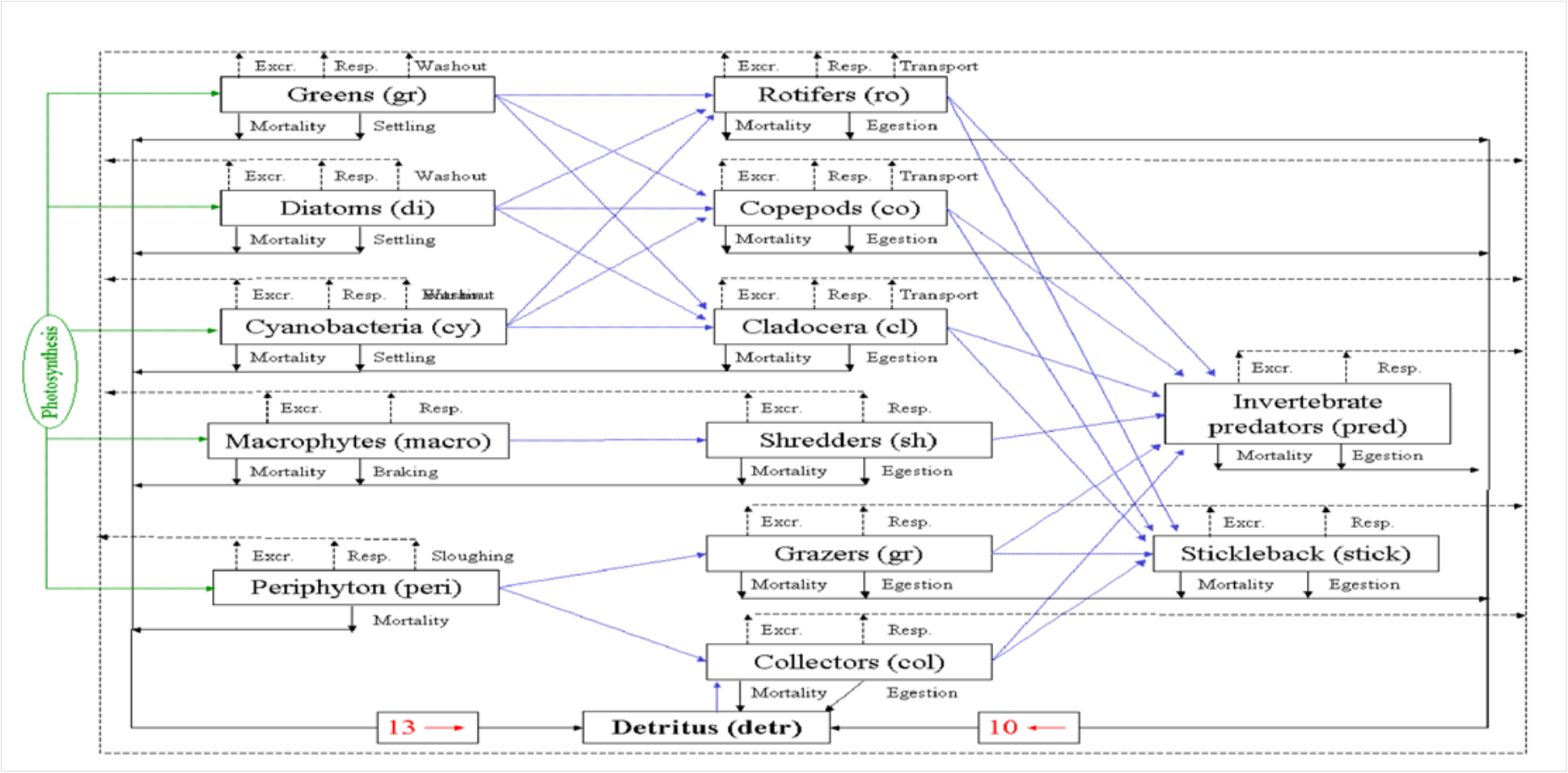}$$
\caption{Representation of the module chain of the aquatic ecosystem
model MELODY.}\label{melody}
\end{figure}

Following this rationale, an aquatic ecosystem model, MELODY\footnote{modelling
MEsocosm structure and functioning for representing LOtic
DYnamic ecosystems}, was built so as to simulate the
functioning of aquatic mesocosms as well as the impact of toxic substances
on the dynamics of their populations. 
A main feature of this kind of ecological model is,
however, the great number of parameters involved in the
modelling: MELODY has a total of
13 compartments  and 219 
parameters; see Figure \ref{melody}. 
These are generally highly uncertain because of both
natural variability and lack of knowledge.
Thus, sensitivity analyses appears an
essential step to identify
non-influential parameters \citep{circif12}. 
These can then be fixed at a nominal value
without significantly impacting the
output.
Consequently, the calibration of the model becomes
less complex \citep{salrat06}.

By a preliminary sensitivity analysis of the
periphyton-grazers submodel (representative of
processes involved in dynamics of primary producers and primary consumers and
involving $20$ input parameters), \cite{ioopop12} concludes that
significant interactions of large degrees (more than three)
exist in this model.
Therefore, a DoCE has to possess excellent subprojection
properties to capture the interaction effects.

\section{Space filling criteria and designs}
\label{sec:designs}

Building a DoCE consists in generating a
matrix $\mathbf{X}^N_d=(x_j^{(i)})_{i=1..N,j=1..d}$,
where $N$ is the number of experiments and $d$ the
number of scalar inputs.
The $i^{th}$ line $x^{(i)}$ of this matrix 
will correspond to the inputs of the $i^{th}$ code execution.
Let us recall that, here,
the purpose is to design $N$ experiments\footnote{
In the following, a ``point'' corresponds to an ``experiment'' (at least a
subset of experimental conditions).}
$x^{(i)}$ to fill as ``homogeneously'' as possible the set $\mathcal{X}={[0;1]}^d$;
even if, in an exploratory phase, a joint probability distribution is not explicitly given
to the inputs, one can consider that this distribution is uniform over
$\mathcal{X}$. In fact, some space-filling criteria discussed hereafter do not
formally refer to the uniform distribution (but geometric considerations), yet they
share with it the property of not favouring any particular area of $\mathcal{X}$.

The most common sampling method is indisputably the
classical Monte Carlo (Simple Random Sampling, SRS), mainly because
of its simplicity and generality \citep{gen03}, but also because of
the difficulty to sample in a more efficient manner when $d$ is large as well.
In our context, it would consist in
randomly, uniformly and independently sampling $d \times N$ draws in $[0;1]$. 
SRS is largely used to propagate uncertainties through computer
codes \citep{derdev08}
because it has the attractive property of giving a convergence rate for
estimates of quantities like expectations that is $O(N^{-1/2})$.
This convergence rate is independent of $d$, the dimension of the input space.

Yet, SRS is known to possess poor space filling properties: SRS leaves wide
unexplored regions and may draw very close points.
Quasi-Monte Carlo sequences \citep{lemi09} furnish good space filling properties,
giving a better convergence rate ($O((\log N)^d/N)$) for estimates of expectations.
If, in practice, improvement over SRS is not attained when $d$ is large,
dimensionality can be controlled by methods of dimension reduction such as
stratified sampling \citep{chedav89} and sensitivity analysis \citep{varlef12}.

The next sections are mainly dedicated to (optimized) LHS, owing to their
property of parsimony mentioned in the introduction, and do not refer to
other interesting classes of SFD, neither maximum entropy designs
nor point process designs in particular. The former is based
on a criterion (Shannon entropy) to maximize which could be used to optimize
LHS.
The resulting SFD is similar to point-distance
optimized LHS (section \ref{sec:point-distance}), as shown by theoretical
works \citep{promul12} and confirmed by our own experiments.
The latter way of sampling seems hardly compatible with LHS
and suffers from the lack of efficient rules
to set its parameters \citep{fra08}.

In this section, criteria used hereafter to make diagnosis of DoCE or to
optimize LHS are defined. Computing optimized LHS requires an efficient,
in fact specialized, optimization algorithm. Since the literature
provides numerous ones, a brief overview of a selection of such algorithms
is proposed. The section ends with our feedback on their behaviours,
with an emphasis on maximin DoCE.

\subsection{Latin Hypercube Sampling}
\label{sec:lhs}

Latin Hypercube Sampling, which is an extension
of stratified sampling, aims at ensuring that each of the
scalar input has the whole of its range
well scanned, according to a probability
distribution\footnote{The range is the support of the distribution.}
\citep{mckbec79,fanli06}.
Below, LHS is introduced in our particular context,
that is considering a uniform distribution over ${[0;1]}^d$, for the sake of
clarity.

\begin{figure}[!ht]
\begin{center}
\begin{picture}(0,110)
\put(-170,135){\includegraphics[scale=0.4,angle=-90]{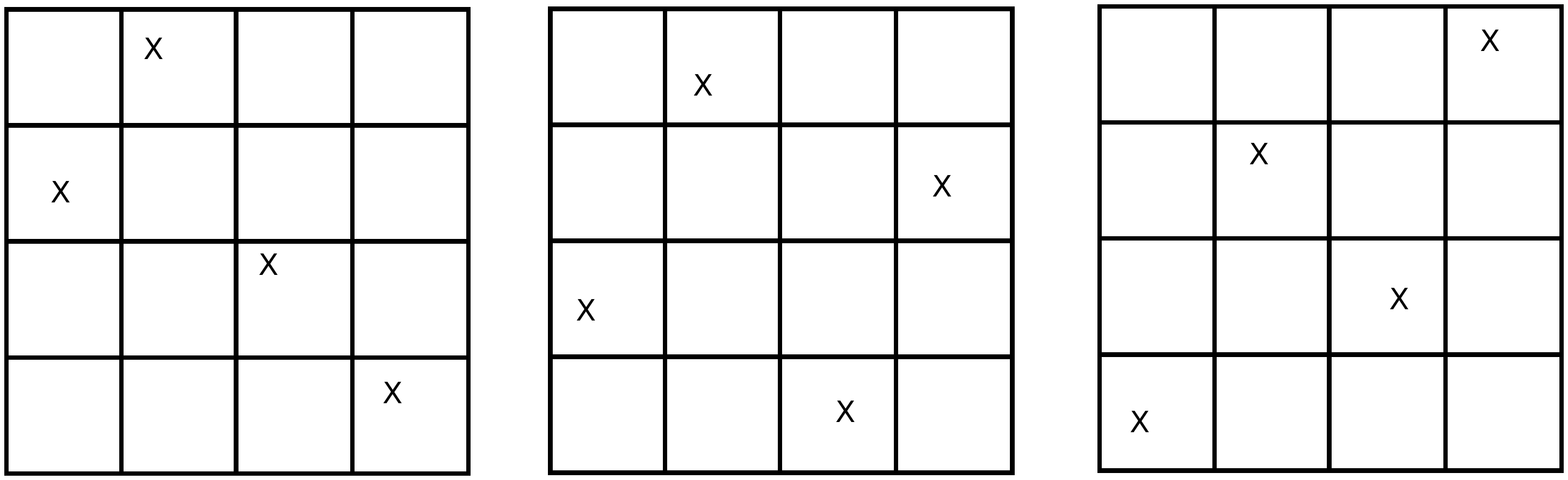}}
\put(-190,50){$X_1$}
\put(-125,-15){$X_2$}
\end{picture}
\caption{Three examples of LHS of size $N=4$
over ${[0;1]}^2$ (with regular intervals): each of the $N$ rows
(each of the $N$ columns, respectively), which corresponds to an interval of
$X_1$ (of $X_2$, resp.), contains one and only one draw $\bx^{(i)}$ (cross).}
\label{fig:LHS}
\end{center}
\end{figure}

Let the range $[0;1]$ of each input variable $X_j$, $j = 1\ldots d$,
be  partitioned into $N$ equally probable intervals
$\mathcal{I}_k = (\frac{k-1}{N}~; \frac{k}{N})$, $k = 1\ldots N$.
The definition of a LHS requires $d$ permutations $\pi_j$ of $\{1,\ldots,N\}$\footnote{
A permutation $\pi$ of $\{1,\ldots,N\}$ is a bijective function from $\{1,\ldots,N\}$ to
$\{1,\ldots,N\}$.}, which are
randomly chosen (uniformly among the $N!$ possible permutations).
The jth component $x_j^{(i)}$ of the ith draw of a (random) LHS is obtained by
randomly picking a value in $\mathcal{I}_{\pi_j(i)}$ (according to the uniform
distribution over $\mathcal{I}_{\pi_j(i)}$). For instance, the second LHS of
Figure~\ref{fig:LHS} (middle) is derived from $\pi_1 = (3; 1; 4; 2)$ and
$\pi_2 = (2; 4; 1; 3)$ (where, by abuse of notation,  each permutation $\pi$
is specified by giving the N-uplet $(\pi(1),\ldots,\pi(N))$).

\subsection{Space filling criteria}
\label{sec:SFC}

As stated previously, LHS is a relevant way to design experiments,
considering one-dimensional projection. Nevertheless, LHS
does not ensure to
fill the input space properly. Some LHS can indeed be really unsatisfactory,
like the first design of Figure~\ref{fig:LHS} which is almost diagonal.
LHS may consequently
perform poorly in metamodel estimation and prediction of the model output \citep{ioobou10}.
Therefore, some authors have proposed to enhance LHS not to only 
fill space in one-dimensional projection, but also in higher dimensions \citep{par94}.
One powerful idea is to adopt some optimality criterion applied to LHS, such as
entropy, discrepancy, minimax and maximin distances, \textit{etc.}
This leads to avoid undesirable situations, such as designs with close points.

The next sections propose some quantitative indicator of space
filling useful i) to optimize LHS or ii) to assess the quality of a design.
Section~\ref{sec:discrepancy} introduces some discrepancy measures which are
relevant for both purposes. Section~\ref{sec:point-distance} and \ref{sec:MST}
introduce some criteria based on the distances between the points of the design.
The former is about the minimax and maximin criteria, which are relevant for
i) but not for ii), and the latter is about a criterion (MST), which
gives an interesting insight of the filling characteristics of a design, but cannot
reasonably be used for i).

\subsubsection{Uniformity criteria}
\label{sec:discrepancy}

Discrepancy measures
consist in judging the uniformity quality of the design. Discrepancy can be seen
as a measure of the gap between the considered
configuration and the uniform one. The star-discrepancy of a design
$\mathbf{X}_d^N=(\bx^{(i)})_{i =1\ldots N}$
over ${[0;1]}^d$ is defined as
\begin{equation}\label{eq:disc}
  D^*(\mathbf{X}_d^N) = \sup_{\by \in {[0;1]}^d} \left|
\frac{1}{N} \, \#\{\bx^{(i)} \in [\b0,\by]\}
- \displaystyle\prod_{j=1}^d y_j
\right|
\end{equation}
where $\by^{\mbox{T}}=(y_1,\ldots,y_d)$ and
$\#\{\bx^{(i)} \in [\b0,\by]\}$ is the number of design points in
$[\b0,\by]=[0,y_1]\times\cdots\times[0,y_d]$.
This discrepancy measure relies on a comparison between the
volume of intervals of form $[\b0,\by]$
and the number of points within these intervals \citep{hic98}.
Definition (\ref{eq:disc}) corresponds to the greater difference
between the value of the CDF of the uniform distribution over ${[0;1]}^d$
(right term) and the value of the empirical CDF of the design (left term).
In practice, the star discrepancy is not computable because of the
\Linf-norm used in formula~(\ref{eq:disc}). 
A preferable alternative is the \Ltwo-norm \citep{fanli06,lemi09}.
For example, the star \Ltwo-discrepancy can be written
as follows:
\begin{equation}\label{eq:discL2}
  D_2^*(\mathbf{X}_d^N) = \left\{ \bigintsss_{{[0;1]}^d} \left[
\frac{1}{N} \, \#\{\bx^{(i)} \in [\b0,\by]\} - \displaystyle\prod_{j=1}^d y_j
\right]^2 d\by \right\}^{\frac{1}{2}} \;.
\end{equation}
A closed-form expression of this discrepancy,
just involving the design point coordinates, can be
obtained \citep{hic98a}.

Discrepancy measures based on \Ltwo-norm are the most popular in
practice because closed-form easy-to-compute expressions are available.
Two of them have remarkable properties: the centered
\Ltwo-discrepancy ($C^2$) and wrap-around \Ltwo-discrepancy ($W^2$)
\citep{jinche05,fan01,fanli06}. Indeed,
Fang has defined seven properties for uniformity measures
including uniformity of subprojections (a particularity of the so-called
modified discrepancies) and invariance by coordinate rotation as well.
Both $C^2$ and $W^2$-discrepancies are defined, by using different
forms of intervals, from the modified \Ltwo-discrepancy.
Their closed-form  expression have been developed in \cite{hic98,hic98a}:
 \begin{itemize}
 \item
the centered \Ltwo-discrepancy
\begin{equation}\label{disccent}
\begin{array}{rl}
\displaystyle C^2(\mathbf{X}_d^N) =  &  \displaystyle \left(\frac{13}{12}\right)^d -\frac{2}{N}\sum_{i=1}^N\prod_{k=1}^d\left(1+\frac{1}{2}|x_k^{(i)}-0.5|-\frac{1}{2}|x_k^{(i)}-0.5|^2\right)\\
& \displaystyle + \frac{1}{N^2}\sum_{i,j=1}^{N}\prod_{k=1}^{d}\left(1+\frac{1}{2}|x_k^{(i)}-0.5|+\frac{1}{2}|x_k^{(j)}-0.5|-\frac{1}{2}|x_k^{(i)}-x_k^{(j)}|\right) \;,
\end{array}
\end{equation}
\item
the wrap-around \Ltwo-discrepancy
\begin{equation}\label{discwrap}
W^2(\mathbf{X}_d^N) = \left(\frac{4}{3}\right)^d +\frac{1}{N^2}\sum_{i,j=1}^N\prod_{k=1}^d\left[\frac{3}{2}-|x_k^{(i)}-x_k^{(j)}|(1-|x_k^{(i)}-x_k^{(j)}|)\right]\;,
\end{equation}
\noindent which allows suppression of bound effects (by wrapping the unit cube for each coordinate).
\end{itemize}
Designs minimizing a discrepancy criterion are called Uniform Designs (UD).

\subsubsection{Point-distance criteria}
\label{sec:point-distance}

\cite{johmoo90} introduced two distance-based criteria. The first idea consists in minimizing the distance
between a point of the input domain and the points
of the design. The corresponding criterion to minimize is called the
minimax criterion $\phi_{mM}$:
\begin{equation}
\phi_{mM}\left(\mathbf{X}_d^N\right)=\max_{\bx \in {[0;1]}^d} \: \min_{i=1\ldots{}N} \:
\|\bx - \bx^{(i)} \|_{\ell^q}\;,
\end{equation}
with $q=2$ (Euclidean distance), typically.
A small value of $\phi_{mM}$ for a design means that there is no point of the
input domain too distant from a point of the design.
This appears important from the point of view of the Gaussian process metamodel, which is
typically based on the assumption of decreasing correlation of outputs with the
distance between the corresponding inputs. However, this criterion needs the
computations of all the distances between every point of the domain and every
point of the design. In practice, an approximation of $\phi_{Mm}$ is obtained
via a fine discretization of the input domain. However, this approach becomes
impracticable for input dimension $d$ larger than three
\citep{promul12}. $\phi_{mM}$ could be also derived from the
Delaunay tessellation that allows to reduce the computational cost
\citep{promul12}, but dimensions $d$ larger than four or five remain an
issue.

A second proposition is to maximize the minimal distance separating two design points.
Let us note $d_{ij}=||\bx^{(i)}-\bx^{(j)}||_{\ell^q}$, with  $q=2$, typically.
The so-called mindist criterion $\phi_{Mm}$ (referring for example to
\verb+mindist()+ routine from DiceDesign R package) is defined as
\begin{equation}
\phi_{Mm}\left(\mathbf{X}_d^N\right)=\min_{\substack{i,j=1\ldots N\\i \ne j}}\,\,d_{ij} \;.
\end{equation}
For a given dimension $d$, a large value of $\phi_{Mm}$ tends to
separate the design points from each other, so allows a better space coverage.

The mindist criterion has been shown to be easily computable but difficult to
optimize. Regularized versions of mindist have been listed in \cite{promul12},
allowing to carry out more efficient optimization in the class of
LHS. In this paper, we use the $\phi_p$ criterion:
\begin{equation}
\phi_{p}\left(\mathbf{X}_d^N\right)=\left[\sum_{\substack{i,j=1\ldots N\\i<j}}\,\,d_{ij}^{-p}\right]^{\frac{1}{p}} \;.
\end{equation}
The following inequality, proved in \cite{promul12}, shows the asymptotic link
between $\phi_{Mm}$ and $\phi_p$. If one defines
$\underline{\xi}^{\star}_{p}$ as the design 
which maximizes $\phi_p$ and $\xi^{\star}$ as the 
one which maximizes $\phi_{Mm}$, then:
\begin{equation}
\label{phip}
1 \ge \frac{\phi_{Mm}(\underline{\xi}^{\star}_{p})}{\phi_{Mm}(\xi^{\star})}\ge{n \choose 2}^{-1/p}.
\end{equation}
Let $\epsilon$ be a threshold, then (\ref{phip}) implies:
\begin{equation}
\frac{\phi_{Mm}(\underline{\xi}^{\star}_{p})}{\phi_{Mm}(\xi^{\star})}\ge 1-\epsilon\,\,\,\, \textrm{for}\,\,\,\, p\gtrapprox2\frac{\ln n}{\epsilon}.
\end{equation}
Hence, when $p$ tends to infinity, minimizing $\phi_p$ is equivalent
to maximizing  $\phi_{Mm}$.
Therefore, a large value of $p$ is required in practice; $p = 50$, as
proposed by \citet{mormit95}, is used in the following.

The commonly so-called maximin design $\mathbf{X}_d^N$ is the one
which maximizes $\phi_{Mm}$, and minimizes the number of pairs of points exactly
separated by the minimal distance.
LHS which are optimized according to the $\phi_p$
or $\phi_{Mm}$ criteria are called maximin LHS in the following.

\subsubsection{Minimum Spanning Tree (MST) criteria}
\label{sec:MST}

The MST \citep{dusras86}, recently introduced
for studying SFD \citep{fravas09,fra08}, enables to analyze the geometrical profile of
designs according to the distances
between points. Regarding design points as vertices, a MST is a tree
which connects all the vertices together and whose sum of edge lengths is minimal.
Once the MST of a design is built,
the mean $m$ and standard deviation $\sigma$ of edge lengths can
be calculated. Designs described as \textit{quasi-periodic}
are characterized by large mean $m$ and small $\sigma$ \citep{fravas09}
compared to random designs or standard LHS.
Such quasi-periodic designs fill the space
efficiently from the point-distance perspective: large $m$ is related to large
interpoint-distance and small $\sigma$ means that the minimal
interpoint-distances are similar.
Moreover, one can introduce a partial order relation for designs based on MST:
a design $D_1$  fills better the space than a design $D_2$ if $m(D_1) > m(D_2)$
and $\sigma(D_1) < \sigma(D_2)$.

\begin{figure}[!ht]
\begin{minipage}{0.5\linewidth}
$$\includegraphics[scale=0.45]{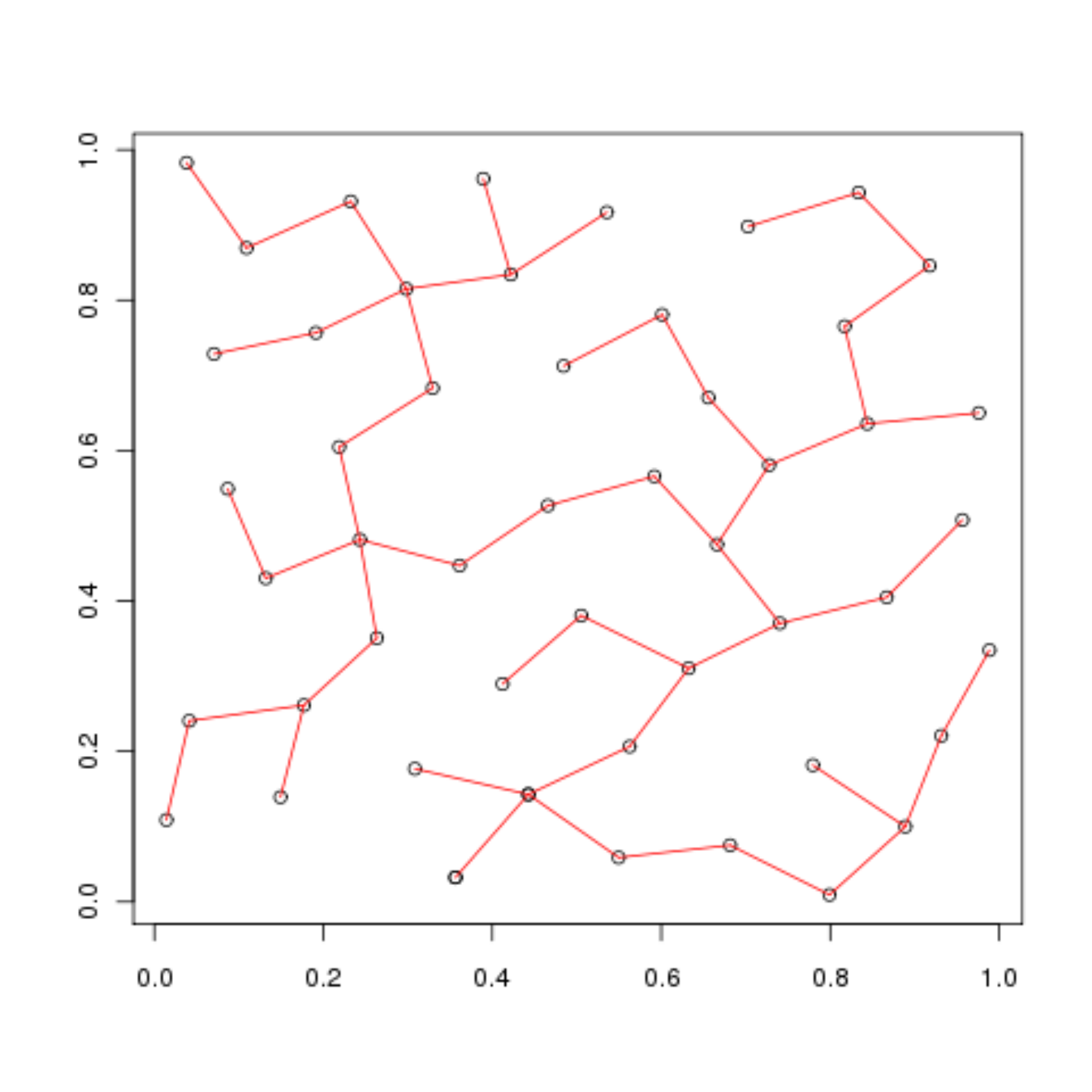}$$
\end{minipage}
\begin{minipage}{0.5\linewidth}
$$\includegraphics[scale=0.60]{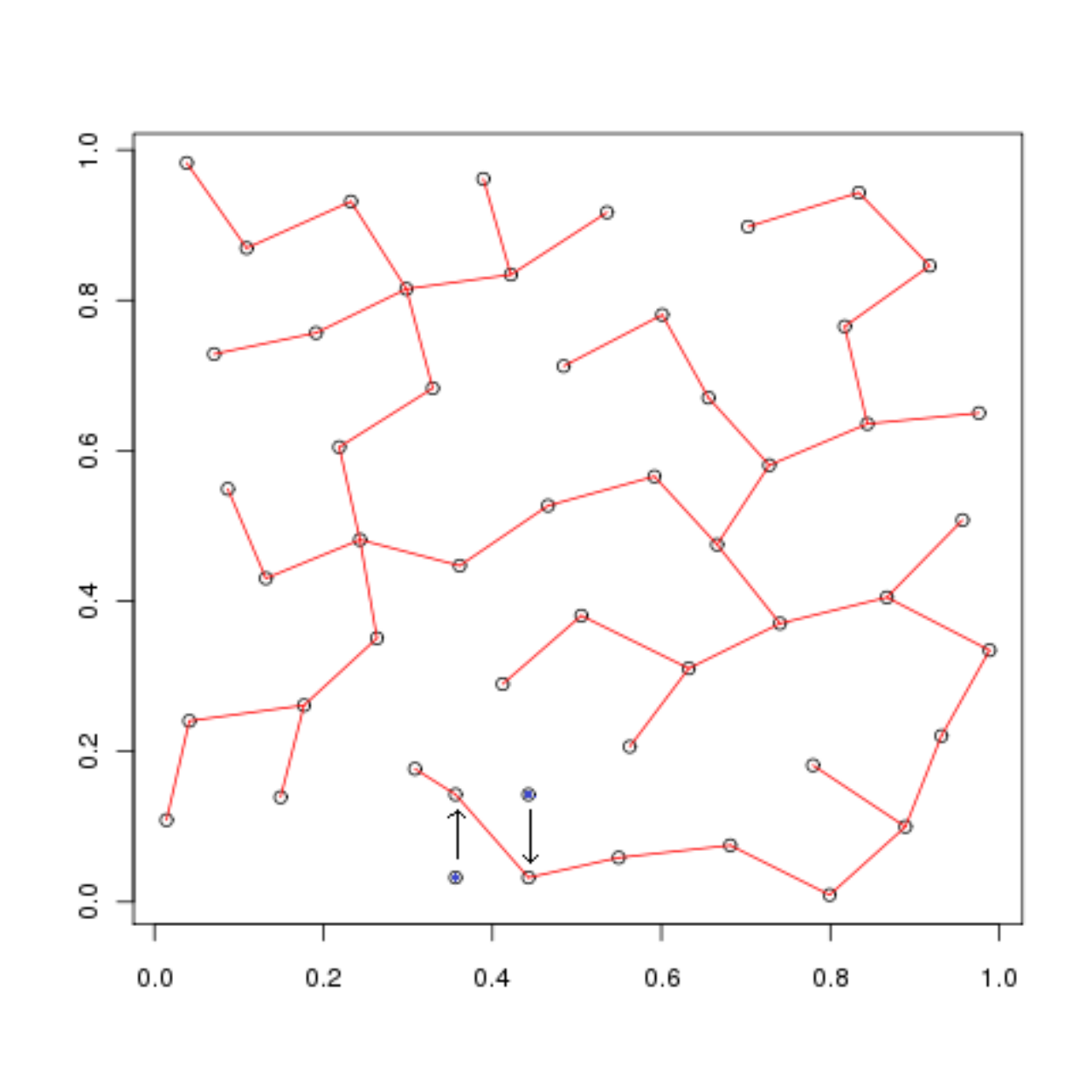}$$
\end{minipage}
\caption{Illustration of two quasi-periodic LHS (left: design with good (large) mindist,
right: design with bad (small) mindist).}
\label{LHSALM}
\end{figure}

MST is a relevant approach focusing on design arrangement and, because $m$ and $\sigma$ are global
characteristics, it makes much more robust diagnosis
than the mindist criterion does. If a design with high
mindist implies a quasi periodic distribution, the reciprocal
is false as illustrated in
Figure~\ref{LHSALM} (on the left). Besides, the MST
criteria appear rather difficult to optimize using stochastic algorithms unlike
the previous ones (see the next section).
However, our numerical experiments lead to conclude that producing
maximin LHS is equivalent to building a quasi-periodic distribution in the LHS design
class.

\subsection{Optimization of Latin Hypercube Sample}

Within the class of latin hypercube arrangements, optimizing a space-filling
criterion in order to avoid undesirable arrangements (such as the diagonal ones
which are the worst cases, see also Figure \ref{fig:LHS}) appear very relevant. 
Optimization can be performed following different approaches, the most natural
being the choice of the best LHS (according to the chosen criterion) among a
large number (e.g. one thousand).
Due to the extremely large number of possible LHS ($(!N)^d$ for discretized LHS
and infinite for randomized LHS), this method is rather inefficient.
Other methods have been developed, based on columnwise-pairwise exchange
algorithms, genetic algorithms, Simulated Annealing (SA), etc.:
see \cite{viaven10} for a review. 
Thus, some practical issues are which algorithm to use to optimize LHS and how
to set its numerical parameters.
Since an exhaustive benchmark of the available methods (with different
parameterizations for the more flexible ones) is hardly possible, we choose to
focus on a limited number of specialized stochastic algorithms:
the Morris and Mitchell (MM) version of SA \citep{mormit95},
a simple variant of MM developed in \cite{mar08} (Boussouf algorithm)
and a stochastic algorithm called ESE (``Enhanced Stochastic Algorithm'',
\cite{jinche05}).
We compare their performance in terms of different space-filling criteria of
the resulting designs.

\subsubsection{Simulated Annealing (SA)}

SA is a probabilistic metaheuristic
to solve global optimization problems.
The approach can provide a good optimizing point in a large search space.
Here, we would like to explore the space of LHS.
In fact, the optimization is carried out from an initial LHS (standard random LHS)
which is (expected to be) improved through elementary random changes.
An elementary change (or elementary perturbation)
of a LHS $\mathbf{X}^N_d$ is done in switching two
randomly chosen coordinates from a randomly chosen column, which keeps the latin
hypercube nature of the sample. The re-evaluation of the criterion after each
elementary change could be very costly (in particular for discrepancy).
Yet, taking into account that only two coordinates are involved in an
elementary change leads to cheap expressions for the re-evaluation.
Formula to re-evaluate in a straightforward way $\phi_p$ and the $C^2$
discrepancy have been established by \cite{jinche05}. We have extended it to
any \Ltwo-discrepancies (including $W^2$ and star \Ltwo-discrepancy).

The main ideas of SA are the following ones.
Designs which do not improve the criterion (bad designs) can be accepted to
avoid getting trapped around a local optimum. At each iteration, an elementary
change of the current design is proposed, then accepted with a probability
which depends on a quantity $T$ called temperature which evolves from an
initial temperature $\mathbf{T_0}$ according to a certain temperature profile.
The temperature decreases with the iterations and less and less bad designs are
accepted. The main issue of SA is to properly set the initial temperature and
the parameters which define the profile to get a good trade-off between a
sufficiently wide exploration of the space and a quick convergence of the
algorithm. 
Finally, a stopping criterion must be specified.
The experiments hereafter are stopped when the maximum number of elementary
changes is reached, which is useful to compare the different algorithms
or to ensure a stop within a chosen duration. However, more sophisticated
criteria could be more relevant to save computations or to carry on
the optimization while the iterations still result in substantial improvements.

The Boussouf SA algorithm has been introduced in \cite{mar08}.
The temperature follows a decreasing geometrical profile:
$T=c^{i}\times T_0$ at the $i$th iteration with $0<c<1$. 
Therefore, the temperature decreases exponentially with the iterations and $c$
must be set very close to $1$ when the dimension $d$ is high.
In this case, SA can sufficiently explore the LHS designs space if enough
iterations are performed and this criterion tends rapidly to a correct
approximation of the optimum.

The MM (Morris and Mitchell) SA algorithm \citep{mormit95} was initially
proposed to generate maximin LHS. It can be used to optimize discrepancy
criteria, or others as well. Contrary to Boussouf SA,
$T$ decreases only if the criterion is not improved
during a row of $I_{max}$ iterations. 
The algorithm stops if the maximum number of elementary changes is reached
or if none of the new designs of the last $I_{max}$ iterations has been accepted.
Morris \& Mitchell proposed some heuristic rules to set 
the different parameters of the algorithm from $N$ and $d$.
We noticed that these rules do not always perform well: some settings can lead
to relatively slow convergences.

\subsubsection{Enhanced Stochastic Evolutionary algorithm (ESE)}

ESE is an efficient and flexible stochastic algorithm for optimizing
LHS \citep{jinche05}. It relies on a precise control of a quantity similar to
the temperature of SA through an exploration step, then an process
of adjustment.
ESE behaves as an acceptance-rejection method like SA. ESE is formed
from an outer loop and an inner loop.
During each of the $M$ iterations of the inner loop,
$J$ new LHS are randomly generated from the current one.
Then, at each iteration of the outer loop, the temperature
is updated via the acceptance ratio; in contrast with SA,
it may increase from an (outer) iteration to the next.
Finally, during $Q$ iterations of the outer loop, 
$J\times M \times Q$ elementary design changes are performed.
The authors expect ESE to behave more efficiently than
SA to optimize LHS.

\subsubsection{Feedback on LHS optimization}

\begin{figure}[!ht]
   \centering
   \rotatebox{90}{\hspace*{0.2\textwidth}\scriptsize{}mindist}
   \includegraphics[height=0.5\textwidth]{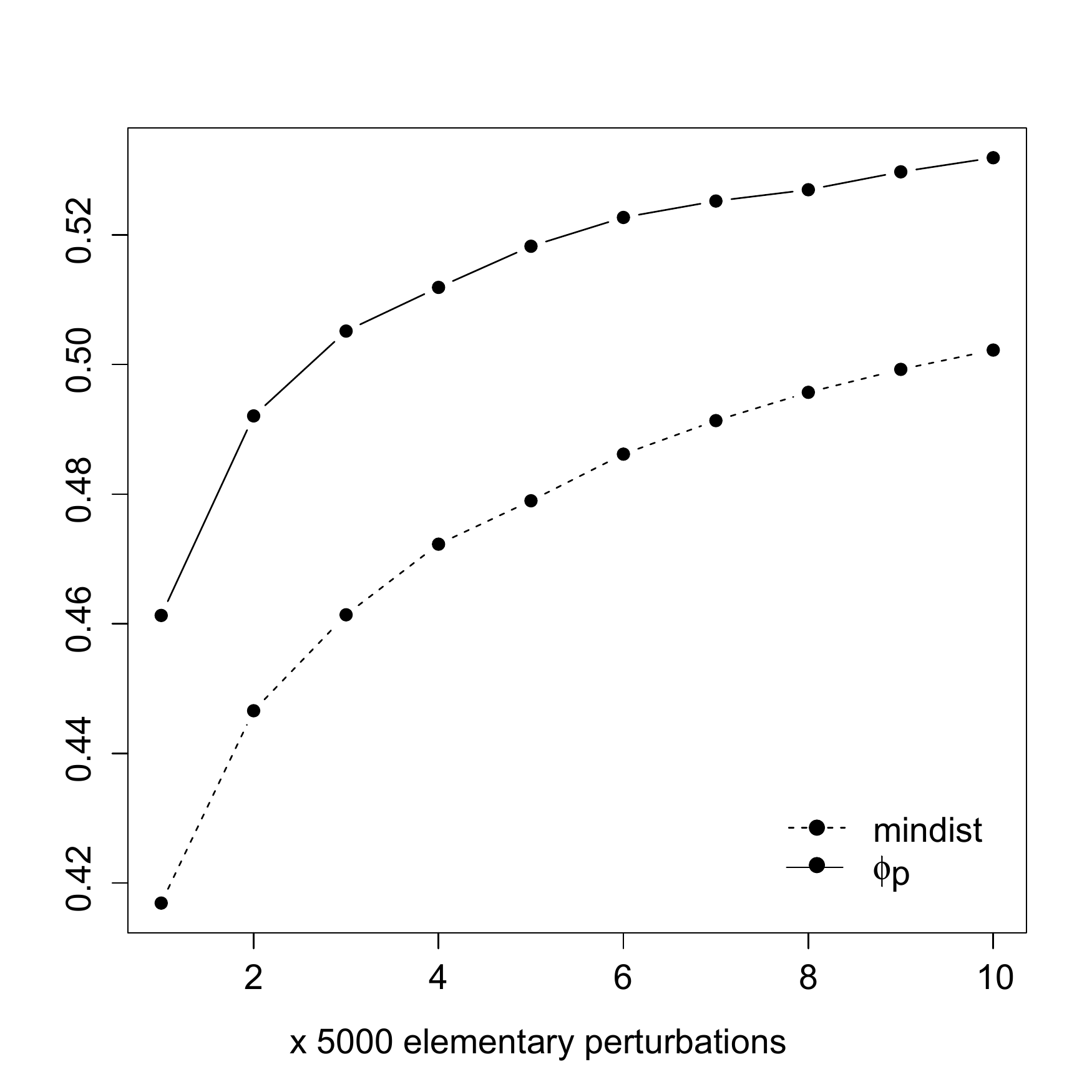}
\caption{Evolution of the mindist criterion in function of the
number of elementary perturbations obtained by ESE ($M=100$, $J=50$)
with $\phi_{p}$ (solid line) and mindist (dashed line);
ordinates are arithmetic means of the best mindist values
from $30$ optimizations of LHS designs ($N=100$, $d=10$).}
\label{lhsOptim}
\end{figure}

First, a test is performed
to show the performance of the regularization of the
mindist criterion (see section \ref{sec:point-distance}).
Using ESE, some optimized LHS are produced by means of
mindist or $\phi_p$; see Figure \ref{lhsOptim}.
If the optimizations are performed with $\phi_p$ instead of mindist,
a clear improvement of mindist is observed.
Hence, the $\phi_p$ criterion is used in the following to build maximin LHS.

Boussouf SA holds a standard geometric profile which is
well-adapted for small values of $d$.
When $d$ rises, making a satisfactory choice of $T_0$ and $c$ gets more and more
difficult. The profile held in MM SA, with $I_{max} \gg 1$
(typically few hundreds),
is actually preferable to perform an efficient exploration of the space.
Thus, the analysis is focused on MM SA
and ESE only in the remainder of this section.

\begin{figure}
\centering
\includegraphics[scale=0.45]{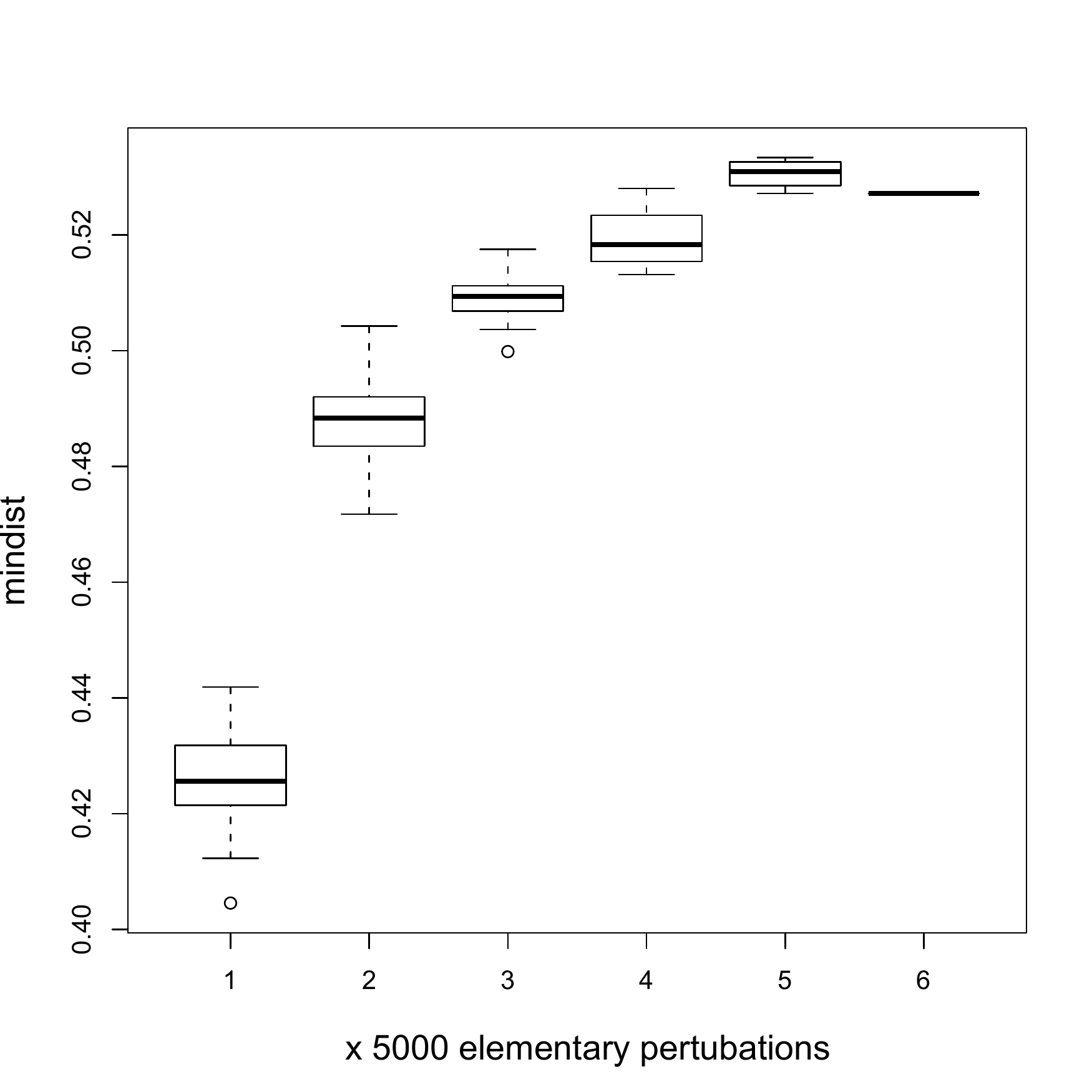}
\includegraphics[scale=0.45]{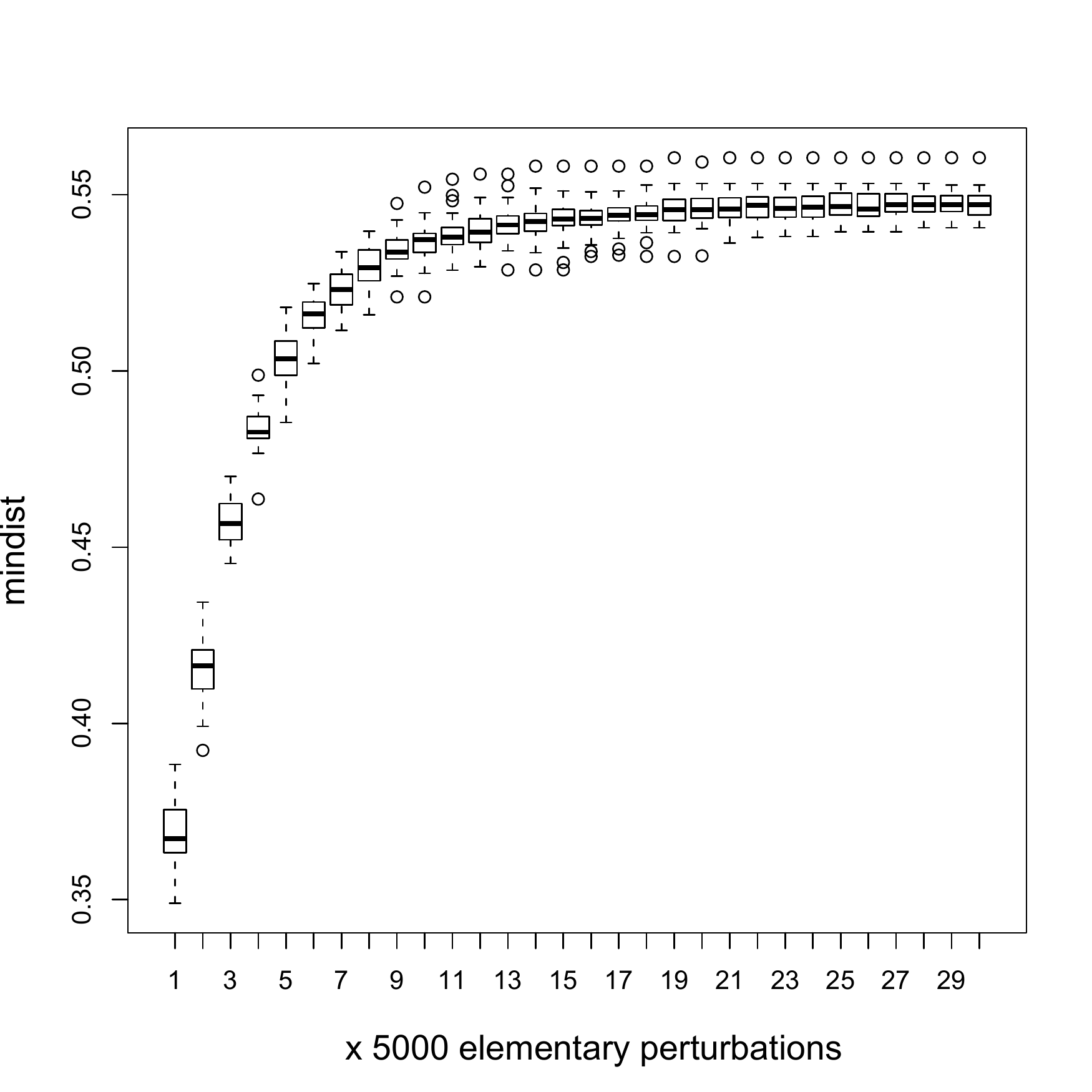}
\caption{Maximin LHS designs ($N=50$, $d=5$) obtained
by MM SA with $\Phi_p$ ($p=50$):
comparison and evolution of the mindist as a function of the number of
elementary changes.
Boxplots are produced from $30$ runs with different initial LHS.
Left (case~$1$): MM SA with $T_0=0.1$, $I_{max}=100$, $c=0.9$. 
Right (case~$2$): MM SA with $T_0=0.1$, $I_{max}=300$, $c=0.9$.}
\label{SAfig}
\end{figure}

The impact of the parameters as a critical aspect affecting
the performances of the algorithms is now considered.
Below, the typical behaviour of MM SA for $d \leq 20$ is illustrated. 
According to our experience, $T_0 = 0.1$ or $0.01$ is a
good choice in order not to accept too many bad designs.
Then, by a care adjustment of
$c$ and $I_{max}$, a trade-off can be found between a relatively
fast convergence and the final optimized value of mindist
(some larger $I_{max}$ and $c$ for a higher dimension $d$).
Figure~\ref{SAfig} shows some results with
two choices of $I_{max}$
to improve the mindist of a $5$-dimensional LHS of size $N=50$:
\begin{itemize}
\item Case $1$ with $T_0=0.1$, $I_{max}=100$ and $c=0.9$ (left of Figure~\ref{SAfig}).
\item Case $2$ with $T_0=0.1$, $I_{max}=300$ and $c=0.9$ (right of Figure~\ref{SAfig}).
\end{itemize}
MM SA converges faster in case~$1$ than in case~$2$ (faster
flattening of the trajectories) and
case~$2$ turns out to be a better choice of $I_{max}$
considering the optimization
of mindist (greater values when convergence seems achieved).

\begin{figure}
\centering
\includegraphics[scale=0.45]{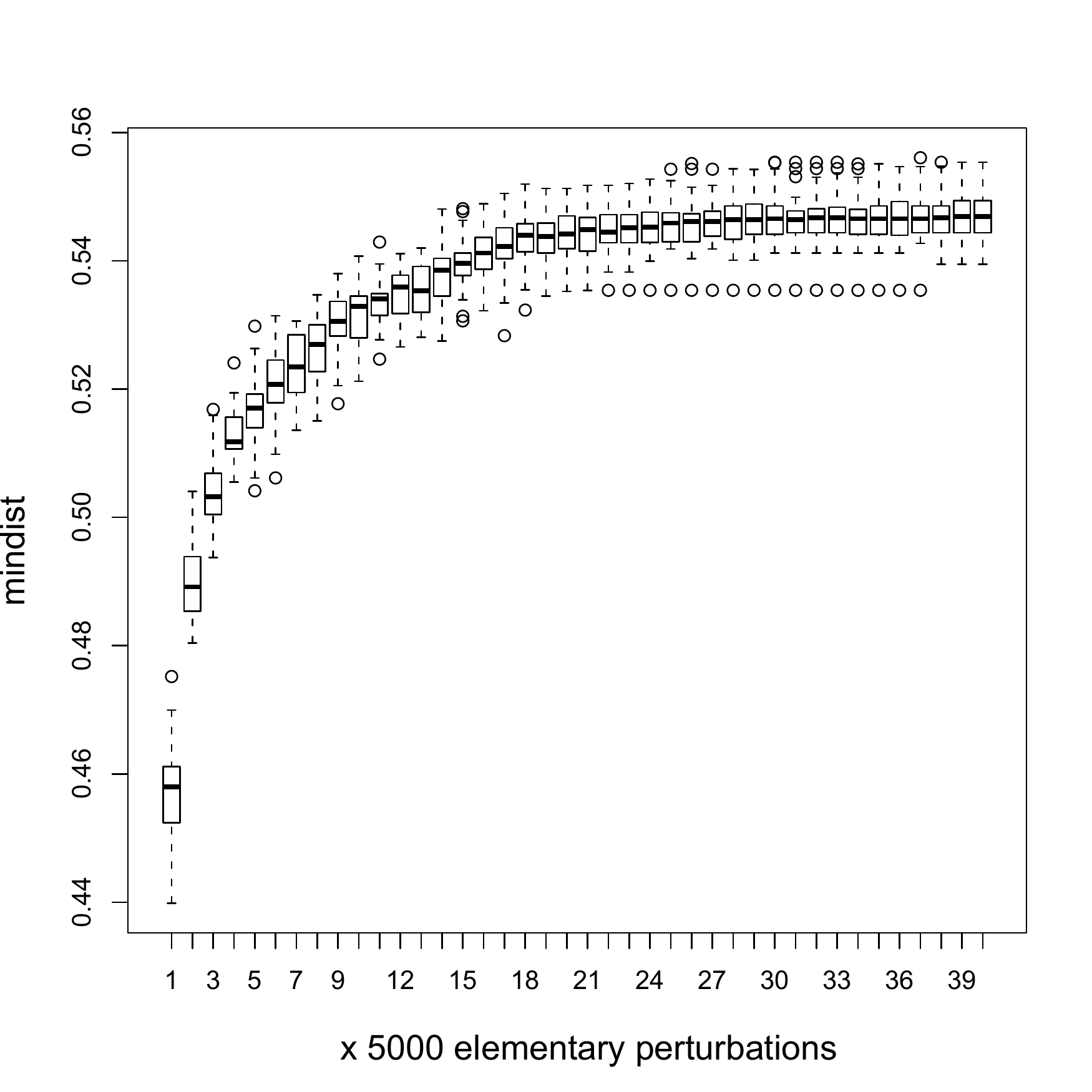}
\caption{Maximin LHS design ($N=50$, $d=5$) obtained by ESE
($M=100$ and $J=50$) with $\Phi_p$ ($p=50$): comparison and evolution
of the mindist criterion
as a function of the number of elementary changes. Boxplots are produced
from $30$ runs with different initial LHS.}
\label{ESEfig}
\end{figure}

ESE is now regarded. In a first attempt, the
parameters suggested in \citet{jinche05} are used: see Figure~\ref{ESEfig}.
Considering results of case~$2$ (right of Figure~\ref{SAfig}),
a faster convergence is obtained with ESE together
with similar final optimized values.
Indeed, most often, only $15,000$ elementary changes are needed
to exceed a mindist of $0.5$ with ESE contrary to what is
observed with MM SA in case~$2$, which is as efficient as
MM SA in case~$1$. Thus, ESE combines the convergence speed
of MM SA in case~$1$ with the ability to find designs of
as large mindist as the outcomes from MM SA in case~$2$.

\begin{figure}[ht!]
\centering
\includegraphics[scale=0.45]{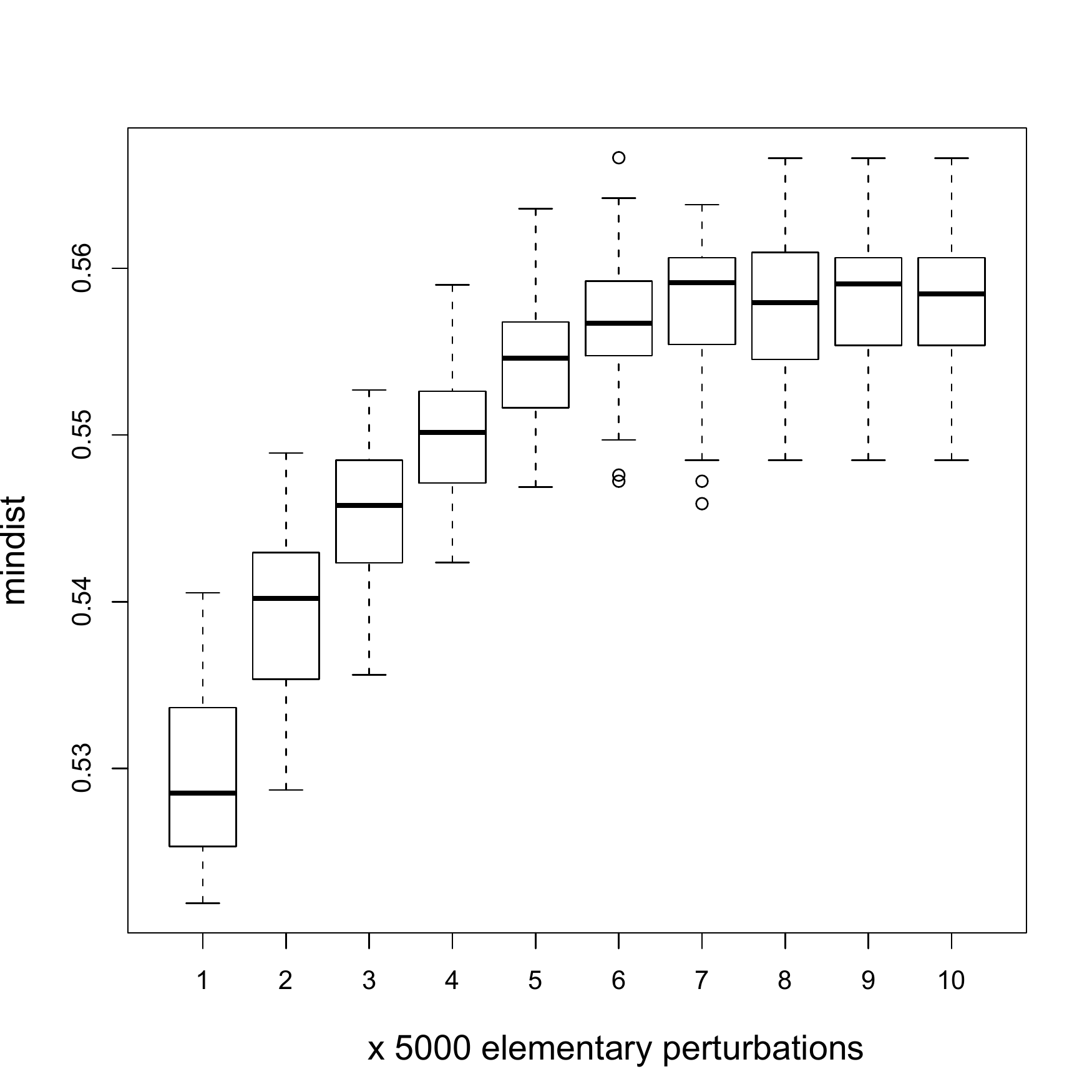}
\includegraphics[scale=0.45]{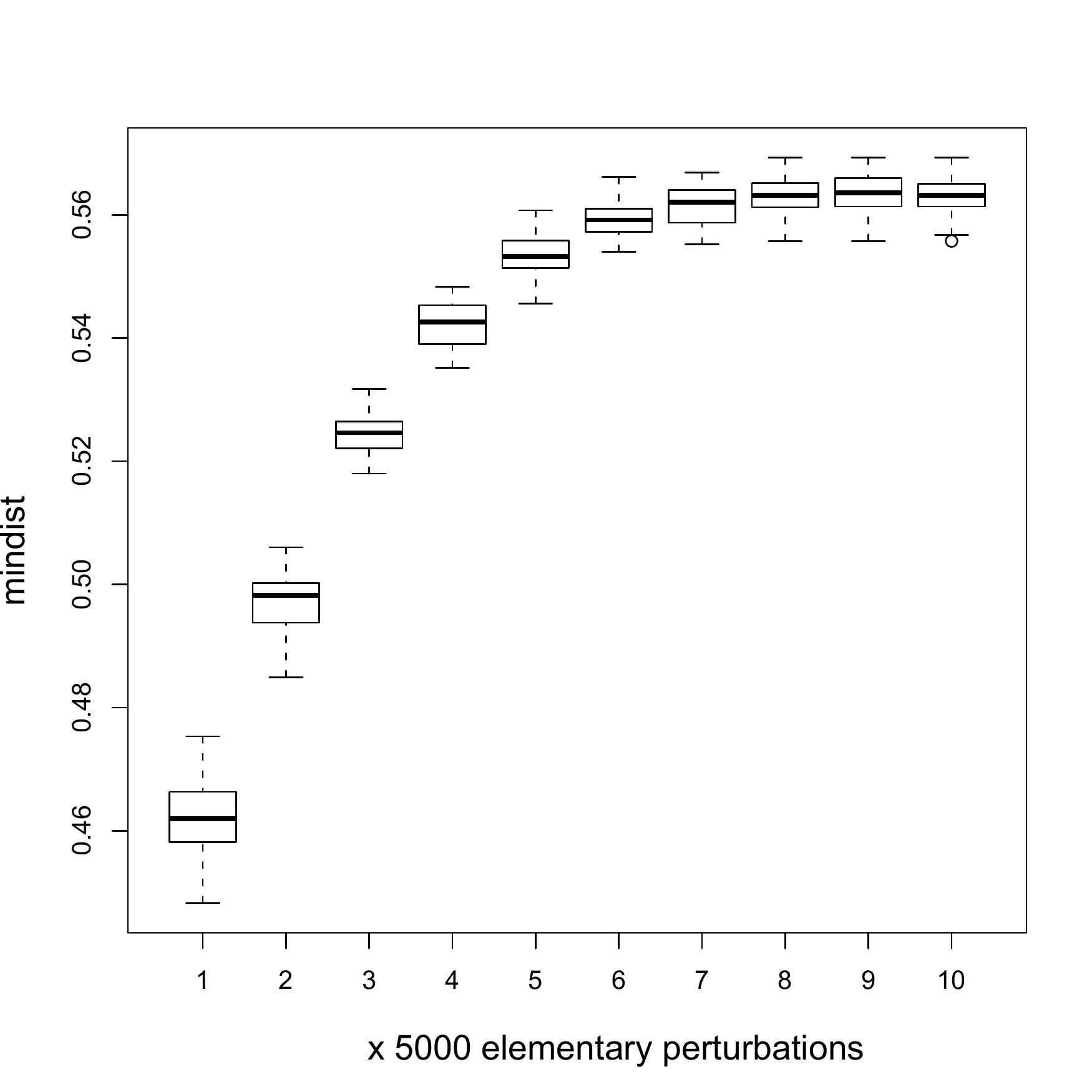}
\caption{Maximin LHS designs ($N=50$, $d=5$) obtained by ESE
and MM SA with $\Phi_p$ ($p=50$): comparison and evolution of the
mindist criterion as a function of the number
of elementary changes. Boxplots are produced through $30$ runs from different
initial LHS. 
Right: MM SA with $T_{0}=0.01$, $I_{max}=1000$ and $c=0.98$.
Left: ESE with $M=300$, $J=50$.}
\label{MMoptfig}
\end{figure}

Finally, the two following tests are proposed:
\begin{itemize}
\item ESE with $M=300$ and $J=50$ (left of Figure~\ref{MMoptfig}).
\item MM SA with $T_0=0.01$, $I_{max}=1000$ and $c=0.98$
(right of Figure~\ref{MMoptfig}).
\end{itemize}
Results of ESE and MM SA are similar.
Both algorithms cost $350,\!000$ elementary
perturbations to exceed $0.56$ which can be considered as a
prohibitive budget in comparison with $15,\!000$ to exceed $0.5$
(Figures~\ref{SAfig} and \ref{ESEfig}).
Generally, the computational budget required to reach an optimized value $v$
of a discrepancy or point-distance based criterion to maximize (respectively,
to minimize) grows much faster than a linear function
with respect to $v$ (resp., of $-v$).

Boussouf SA is used in the following section to carry out comparisons of LHS
optimizations. This choice enables us to perform many optimizations up to relatively
large dimensions within a reasonable duration.
If the results cannot be claimed to be perfectly fair due
to unfinished optimizations, they are representative of the practice of LHS
optimization.

\section{Robustness to projections over 2D subspaces}\label{sec:subproj}
\label{sec:robustness}

\subsection{Motivations}

An important characteristic of a SFD $\mathbf{X}_d^N$ over ${[0;1]}^d$ is
its robustness to projections over lower-dimensional subspaces, which means
that the $k$-dimensional subsamples of the SFD, $k<d$, obtained by deleting
$d-k$ columns of the matrix $\mathbf{X}_d^N$, fill ${[0;1]}^k$
efficiently (according to a space filling criterion).
A LHS structure for the SFD is not sufficient because it only guarantees good
repartitions for one-dimensional projections, and not for projections to
subspaces of greater dimensions.
Indeed, to capture precisely an interaction effect between some inputs, a good
coverage is required onto the subspace spanned by these inputs (see
section~\ref{sec:prey}).

Another reason why that property of robustness really matters is that a
metamodel fitting can be made in a smaller dimension than $d$ (see an example
in \cite{cangar08}).
In practice, this is often the case because the output analysis of an initial
design (``screening step'') may reveal some useless (i.e. non influent) input
variables that can be neglected during the metamodel fitting step \citep{puj09}.
Moreover, when a selection of input variables is made during the metamodel
fitting step (as for example in \cite{marioo07}), the new sample, solely
including the retained input variables, has to keep good space filling
properties.

The remainder of the paper is focused on the space filling quality of the
2D projections of optimized LHS.
Most of the time, interaction effects of order two (i.e. between two inputs)
are significantly larger than interaction of order three, and so on.

\begin{figure}
\begin{center}
\includegraphics[scale=0.45]{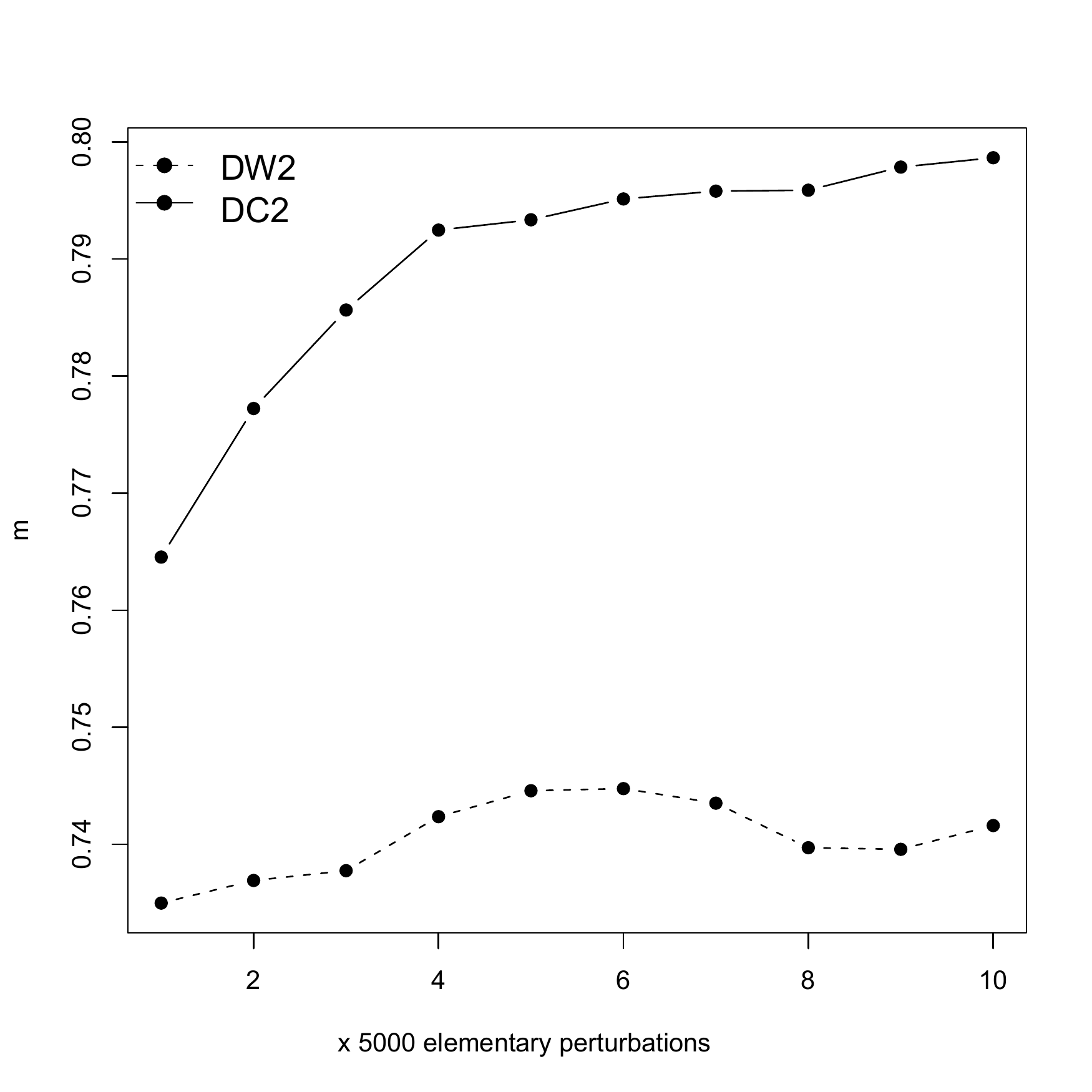}
  \includegraphics[scale=0.45]{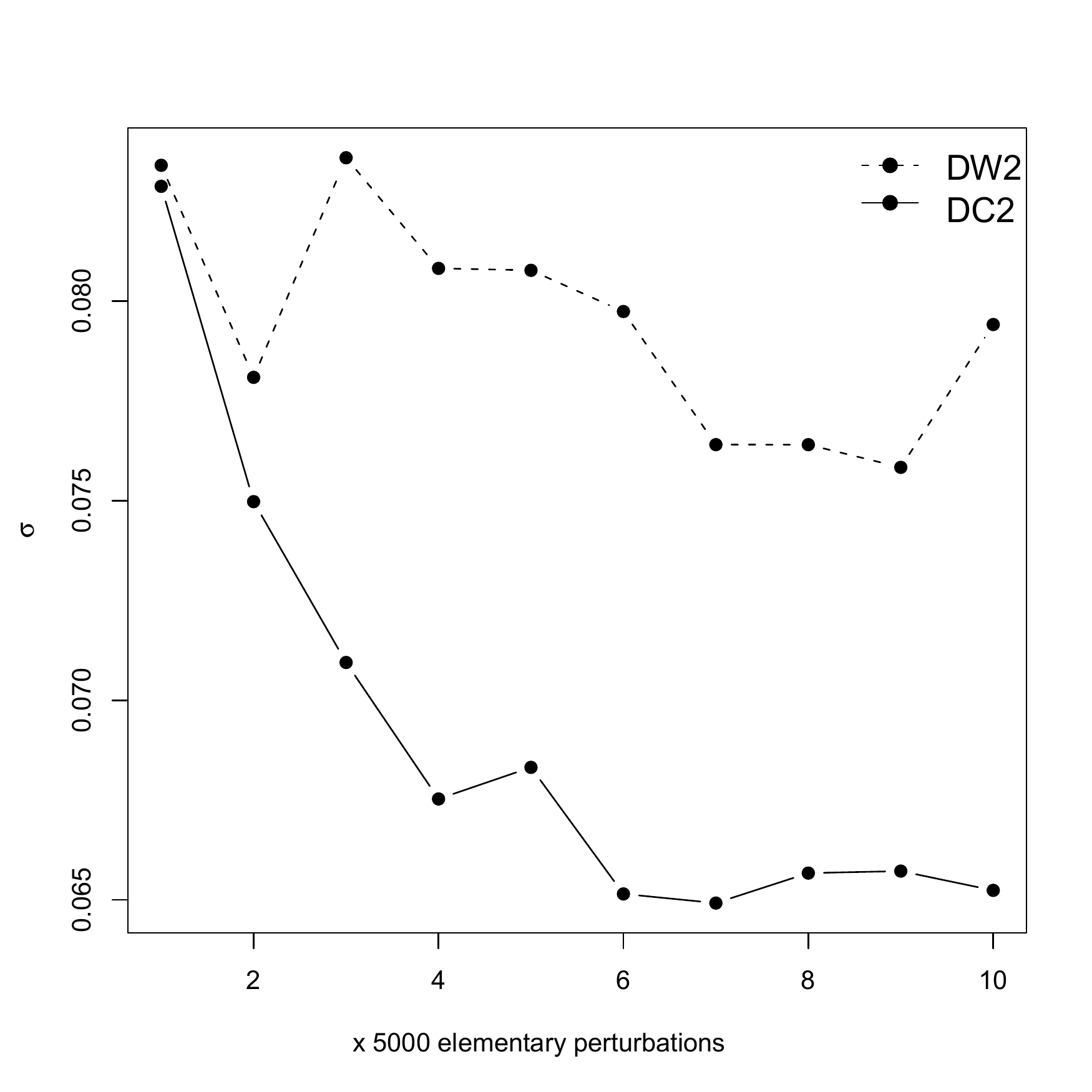}
  \end{center}
\vspace*{-1cm}
\caption{MST mean $m$ (left) and standard deviation $\sigma$ (right)
of $C^2$-discrepancy (solid line) and $W^2$-discrepancy
(dashed line) optimized LHS of size $N=100$ and dimension $d=10$.
Ordinates are arithmetic means of five optimization outcomes.}
\label{MSTC2W2}
\end{figure}

Discrepancy and point-distance based criteria can be regarded as relevant
measures to quantify the quality of space-filling designs. Unfortunately,
they appear incompatible in high dimension in the sense
that a quasi-periodic design (see section~\ref{sec:MST}) does not reach the
lowest values of  any discrepancy measure (for a LHS of same size $N$).  
In high dimension, if a discrepancy-optimized LHS is compared to a
maximin LHS, some differences between their MST criteria, mean $m$ and
standard deviation $\sigma$ (see section~\ref{sec:MST}), are observed. 
Moreover, Figure \ref{MSTC2W2} indicates that means $m$, respectively
standard deviations $\sigma$, of $C^2$-discrepancy
optimized LHS are larger, resp. smaller, than ones of
$W^2$-discrepancy optimized LHS.
This is the reason why we will mainly focus on the
$C^2$-discrepancy rather than the $W^2$-discrepancy
(in the view of the partial order relation defined in
section \ref{sec:MST}).

Below, some tests underline the space filling properties of LHS
2D subprojections.
Qualities of LHS optimized according to different criteria
are analyzed by means of some discrepancies and the MST criteria.
All the tests of sections~\ref{sec:discrepancyAnalysis}
and \ref{sec:MSTanalysis} are performed with
designs of size $N=100$ and of dimension $d$ ranging from $2$ to $54$.
To our knowledge, this is the first numerical study with such a
range of design dimensions.
To capture any variability due to the optimization process, some boxplots
are built from all 2D subsamples of five optimized LHS per dimension.

\subsection{Analysis according to \Ltwo-discrepancy measures}
\label{sec:discrepancyAnalysis}

\begin{figure}[!ht]
\begin{center}
   \includegraphics[scale=0.6]{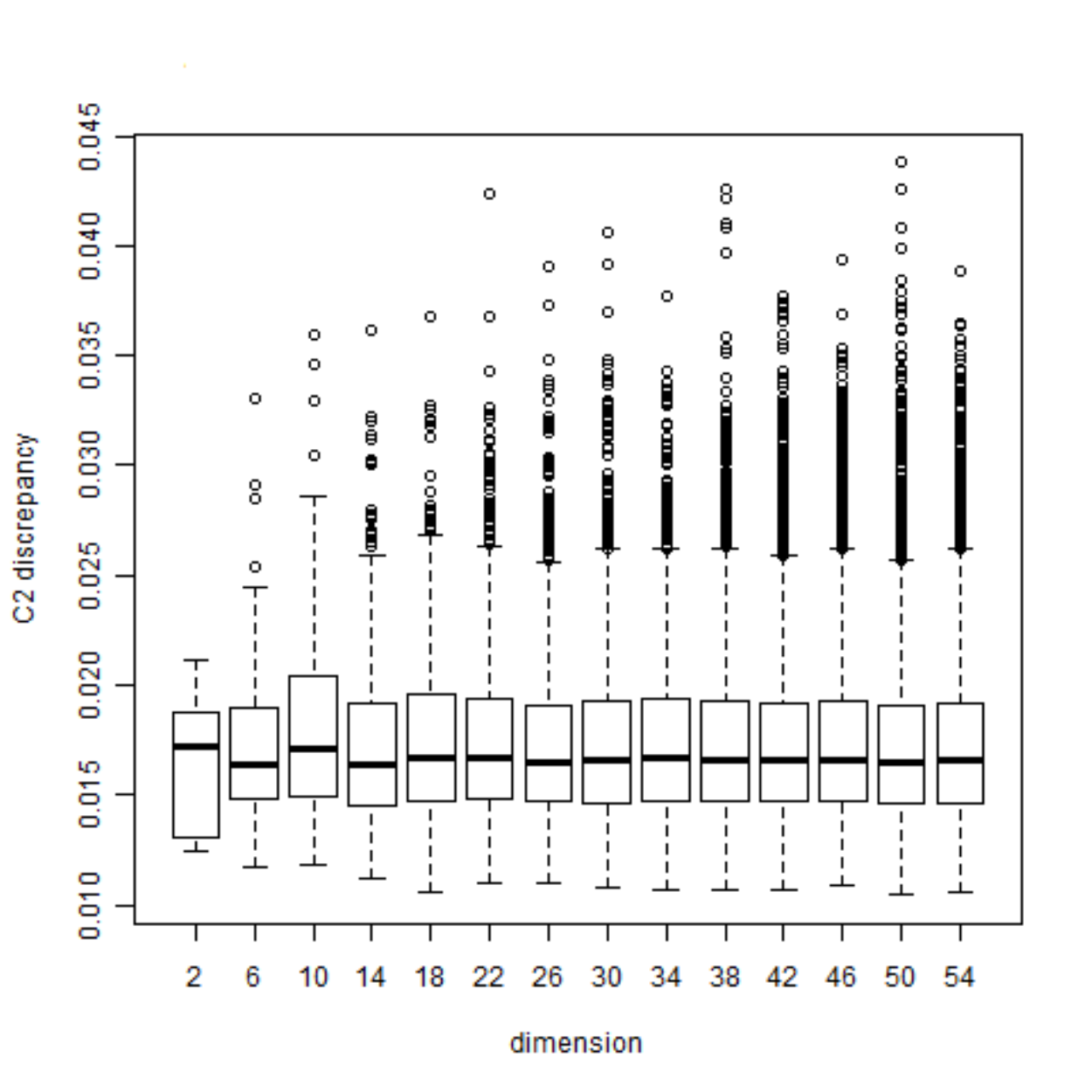}
   \includegraphics[scale=0.6]{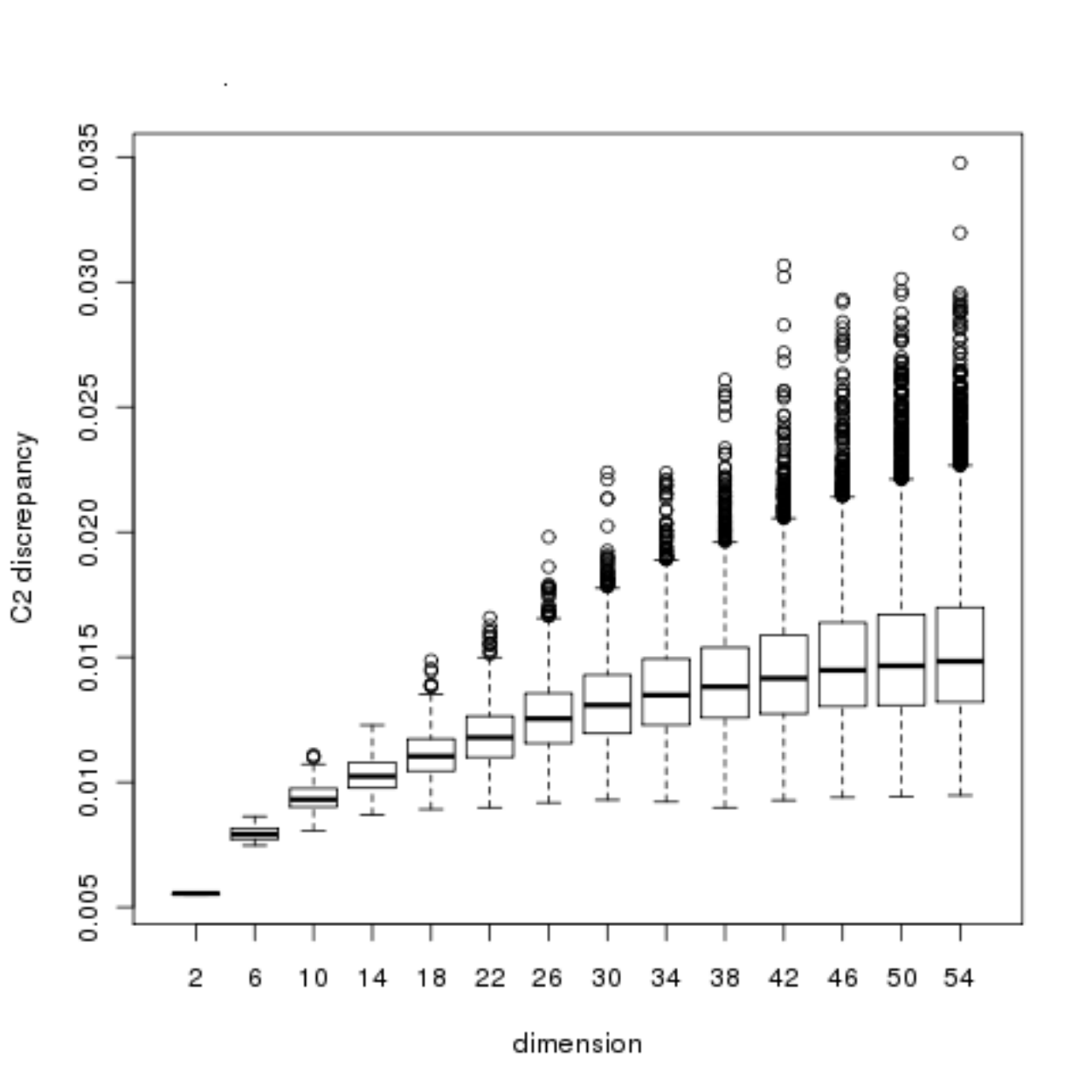}
\end{center}
\vspace*{-1cm}
\caption{$C^2$-discrepancies of the 2D subsamples of
non-optimized random LHS (left) and $C^2$-discrepancy optimized LHS (right).}
\label{projectionsL2}
\end{figure}

\begin{figure}[!ht]
\begin{center}
   \includegraphics[scale=0.6]{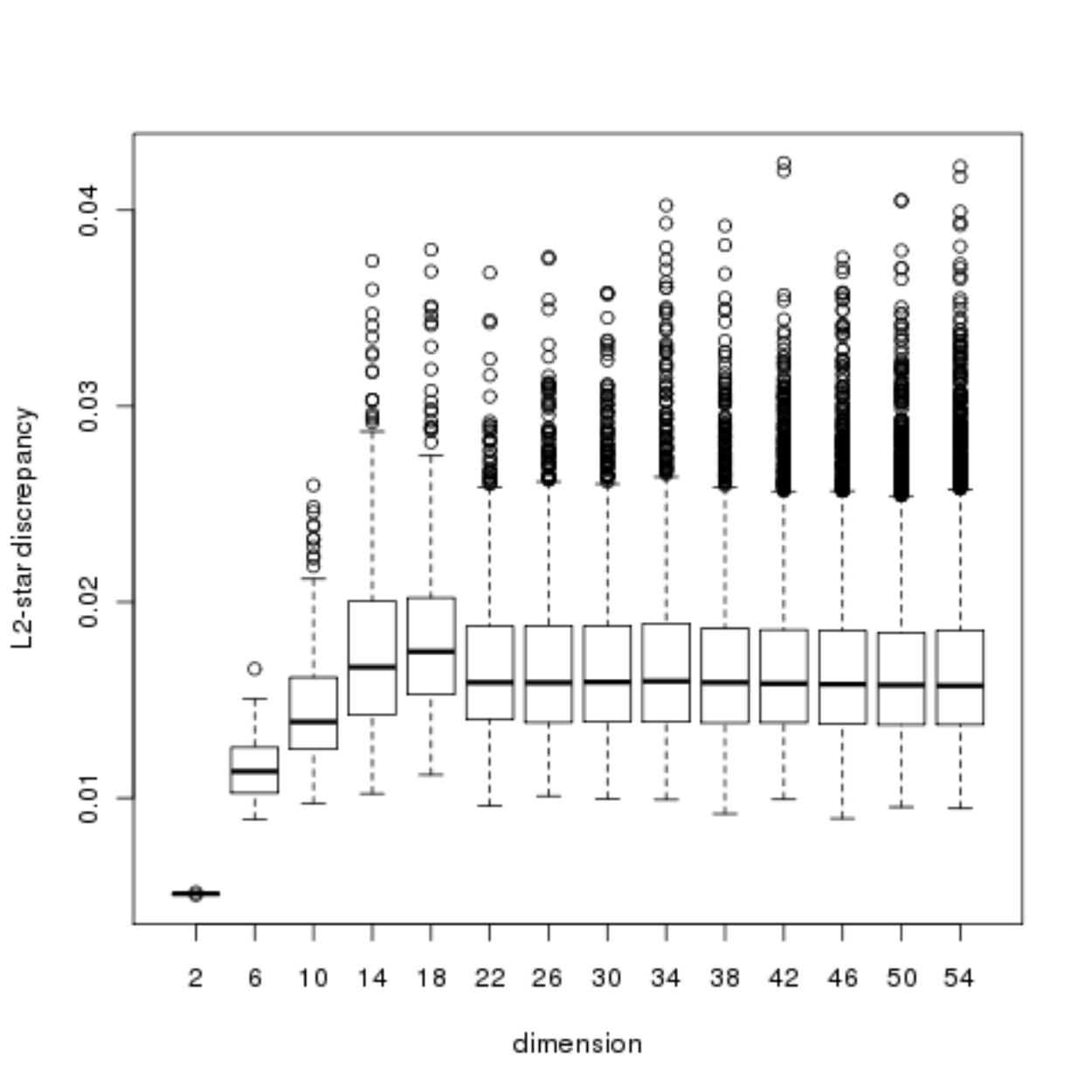}
   \includegraphics[scale=0.6]{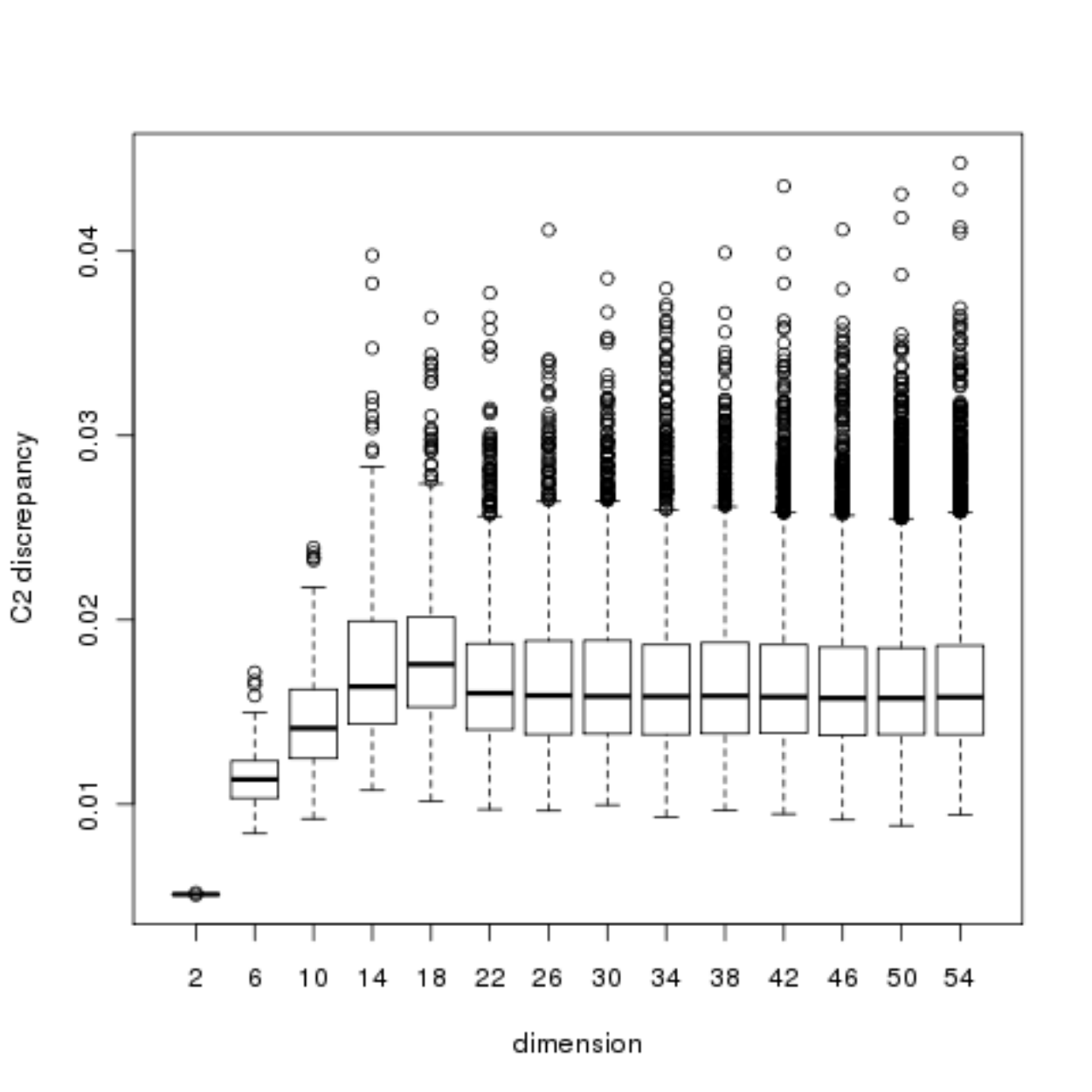}
\end{center}
\vspace*{-1cm}
\caption{\Ltwo-star-discrepancies (left) and $C^2$-discrepancies (right) 
of the 2D subsamples of \Ltwo-star-discrepancy optimized LHS.}
\label{projectionsL2star}
\end{figure}

Figure~\ref{projectionsL2} displays the $C^2$-discrepancies
of non-optimized LHS and $C^2$-optimized LHS.
The comparison of medians
points out that the $C^2$-optimized LHS are robust to 2D subprojections
in the sense that the 2D subsamples get a reduced $C^2$-discrepancy.
The optimizations with the $C^2$-discrepancy appear efficient
for dimensions as large as $54$
since the $C^2$-discrepancies of the 2D subsamples are then less than the
reference value corresponding to the 2D subsamples of the non-optimized LHS
(around $0.017$).
This is fully consistent with Fang's requirement for
the so-called modified discrepancies (see section~\ref{sec:discrepancy}).
Moreover, the increase of $C^2$-discrepancy (median) of the $C^2$-discrepancy
optimized LHS with the dimension $d$ is constant and slow.
Similar observations can be made with the $W^2$-discrepancy
instead of the $C^2$-discrepancy.
Besides, when the same experiment is carried out
with the \Ltwo-star-discrepancy (which is not a modified discrepancy),
results are different: see
Figure~\ref{projectionsL2star}.
In this case, the increase of discrepancy
with $d$ is rather fast and the influence of the optimization
on 2D subsamples ceases for $d > 10$. 

\begin{figure}[!ht]
\begin{center}
\includegraphics[scale=0.6]{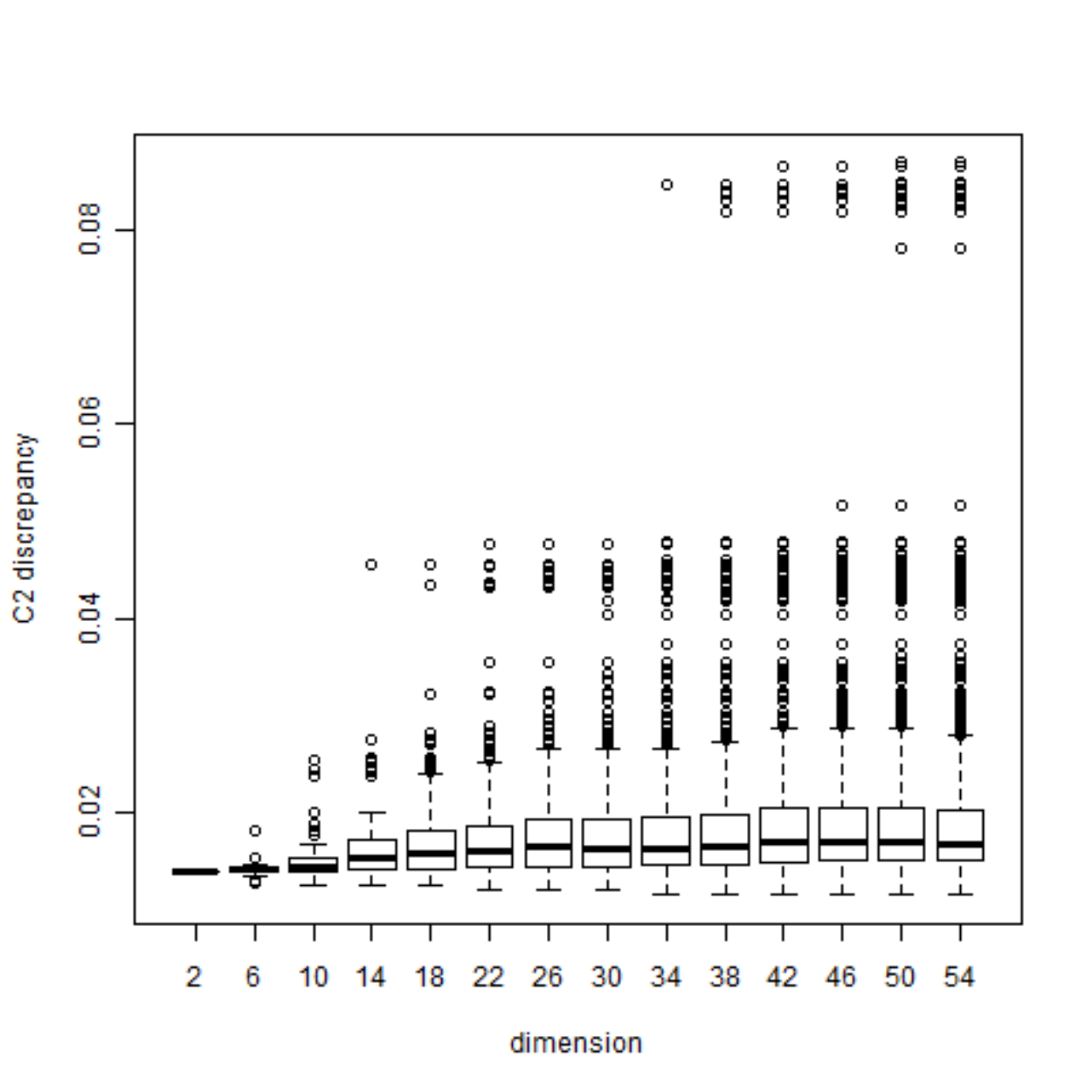}
\end{center}
\vspace*{-1cm}
\caption{$C^2$-discrepancies of the 2D subsamples of the scrambled Sobol' sequence.}
\label{projectionsSobolLHS}
\end{figure}

Figure~\ref{projectionsSobolLHS} illustrates the $C^2$-discrepancies
of the 2D subsamples of the classical low-discrepancy sequence of Sobol,
with Owen scrambling
in order to break the alignments created by this deterministic sequence.
It confirms a well-known fact: the Sobol' sequence has several
poor 2D subprojections in terms of discrepancy, especially
in high dimension.

\begin{figure}[!ht]
\begin{center}
\includegraphics[scale=0.6]{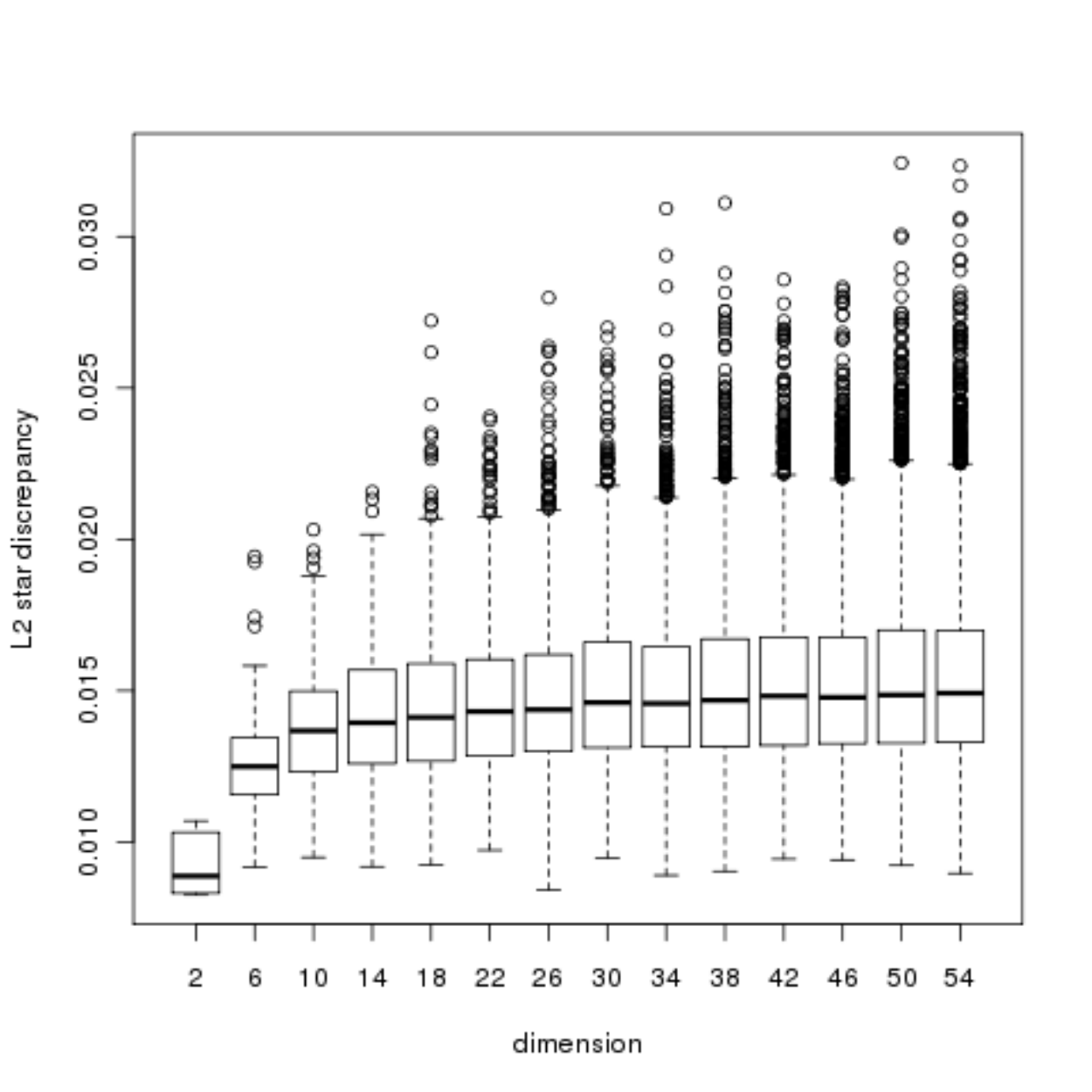}
\includegraphics[scale=0.6]{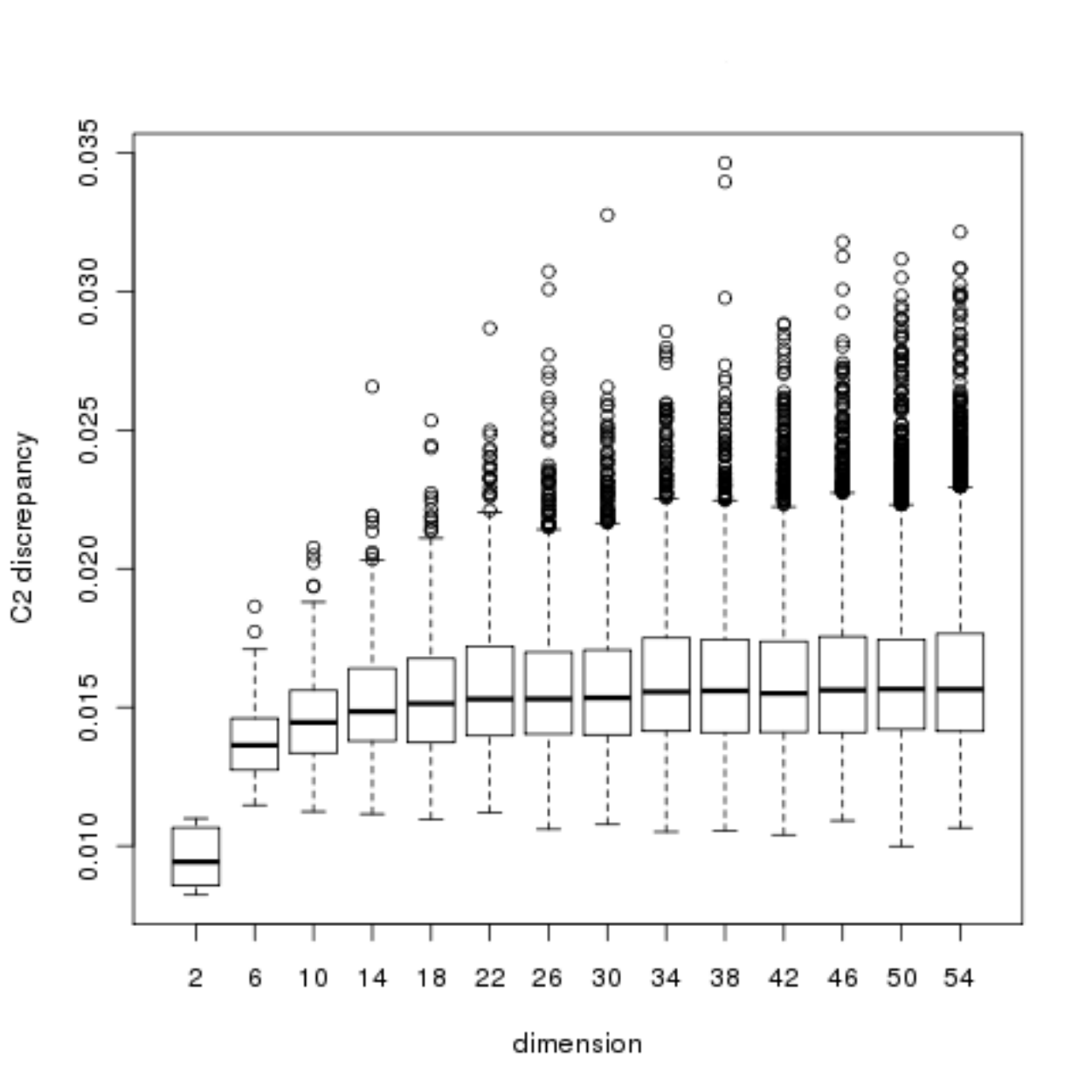}
\end{center}
\vspace*{-1cm}
\caption{
\Ltwo-star-discrepancies (left) and $C^2$-discrepancies (right) 
of the 2D subsamples of maxi\-min LHS.}
\label{projectionsL2b}
\end{figure}

Finally, the same tests are performed on maximin designs.
Previous works \citep{mar08,ioobou10} have shown that maximin
designs are not robust in terms of the mindist criterion
of their 2D subsamples.
It can be explained by alignments in the overall space ${[0;1]}^d$.
Figure~\ref{projectionsL2b} reveals that the maximin design are
not robust to projections to 2D subspaces according to two
different discrepancies, since results are similar to the ones of
the \Ltwo-star-discrepancy optimized LHS (Figure~\ref{projectionsL2star}).
It strongly confirms previous conclusions: the maximin design is not
recommended if one of the objectives is to have some good space
fillings of the 2D subspaces.

\subsection{Analysis according to the MST criteria}
\label{sec:MSTanalysis}

Due to the lack of robustness of the mindist criterion mentioned in
section~\ref{sec:MST}, only the MST criteria are regarded in the
point-distance perspective. 

\begin{figure}
\begin{center}
\includegraphics[scale=0.45]{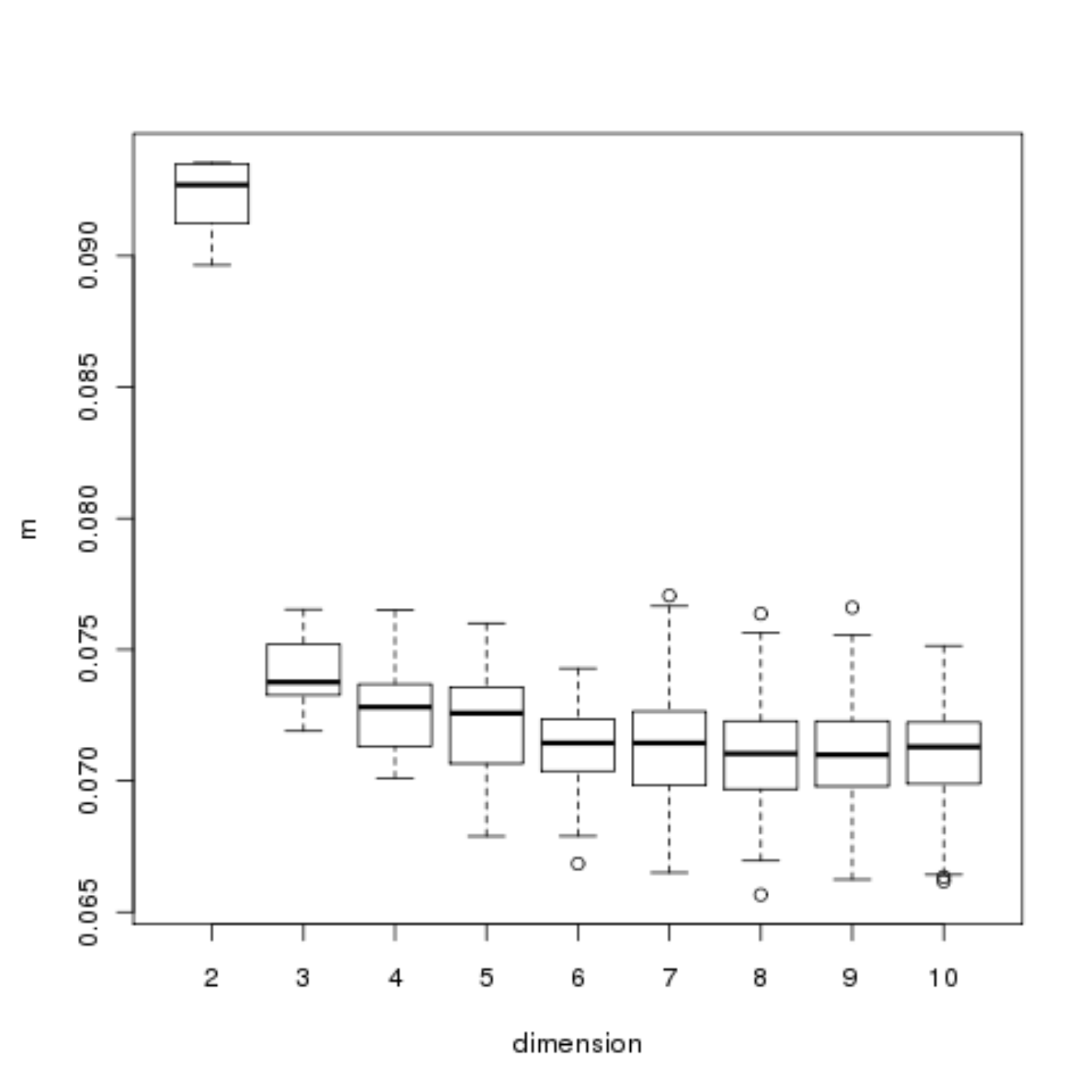}
\includegraphics[scale=0.45]{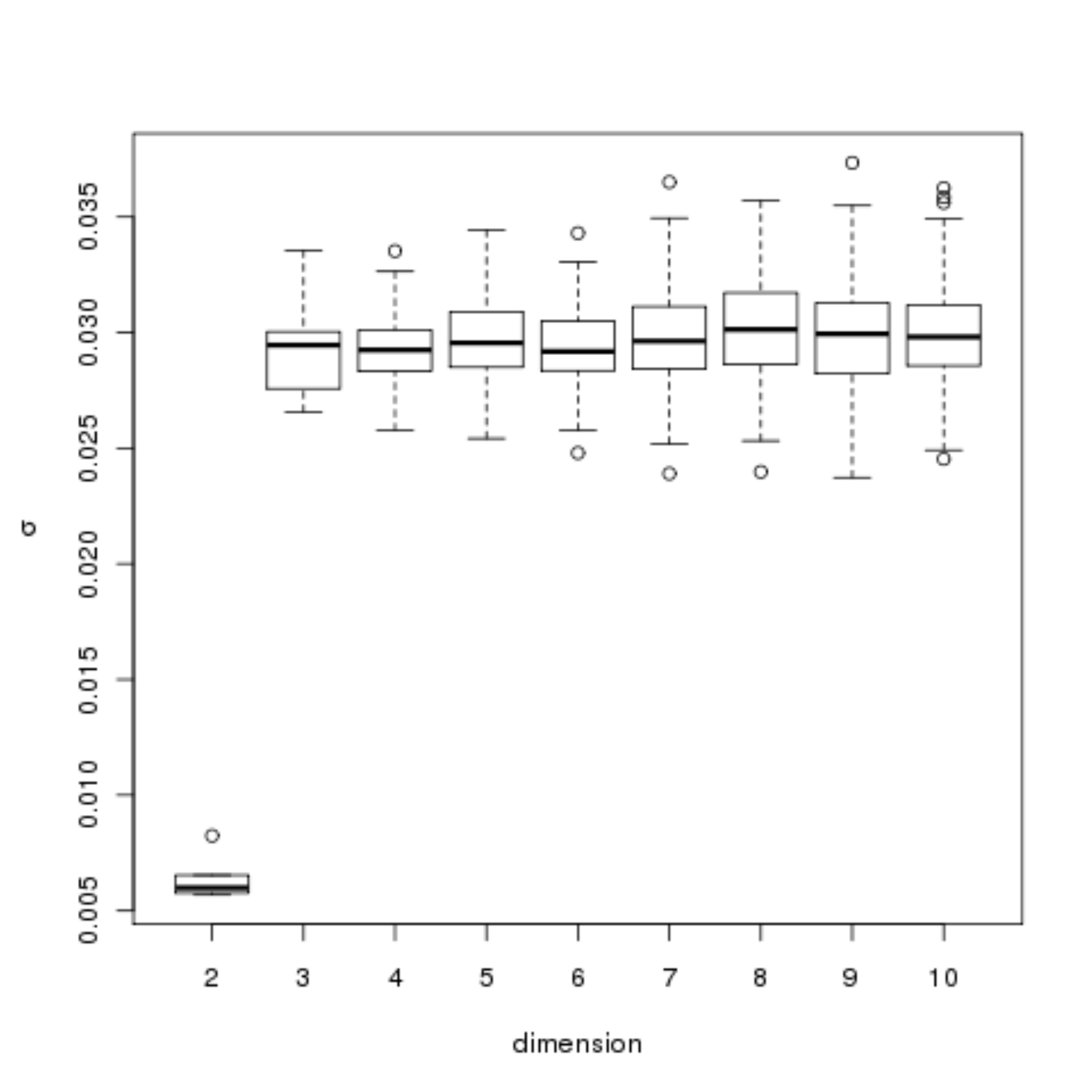}
\end{center}
\vspace*{-1cm}
\caption{MST mean $m$ (left) and standard deviation $\sigma$ (right)
of the 2D subsamples of maximin LHS.}
\label{MSTmaximin}
\end{figure}

\begin{figure}[!ht]
\begin{center}
\includegraphics[scale=0.45]{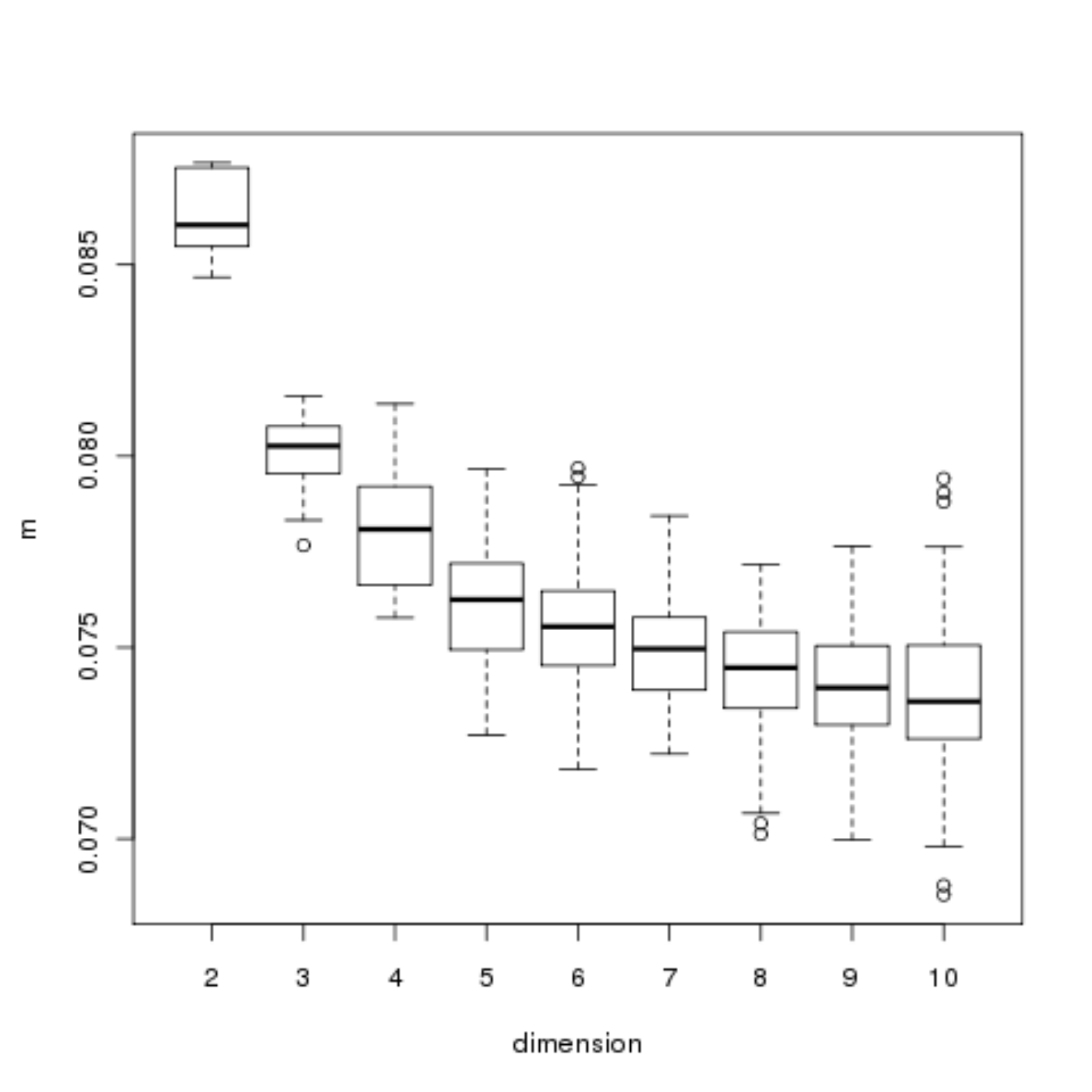}
\includegraphics[scale=0.45]{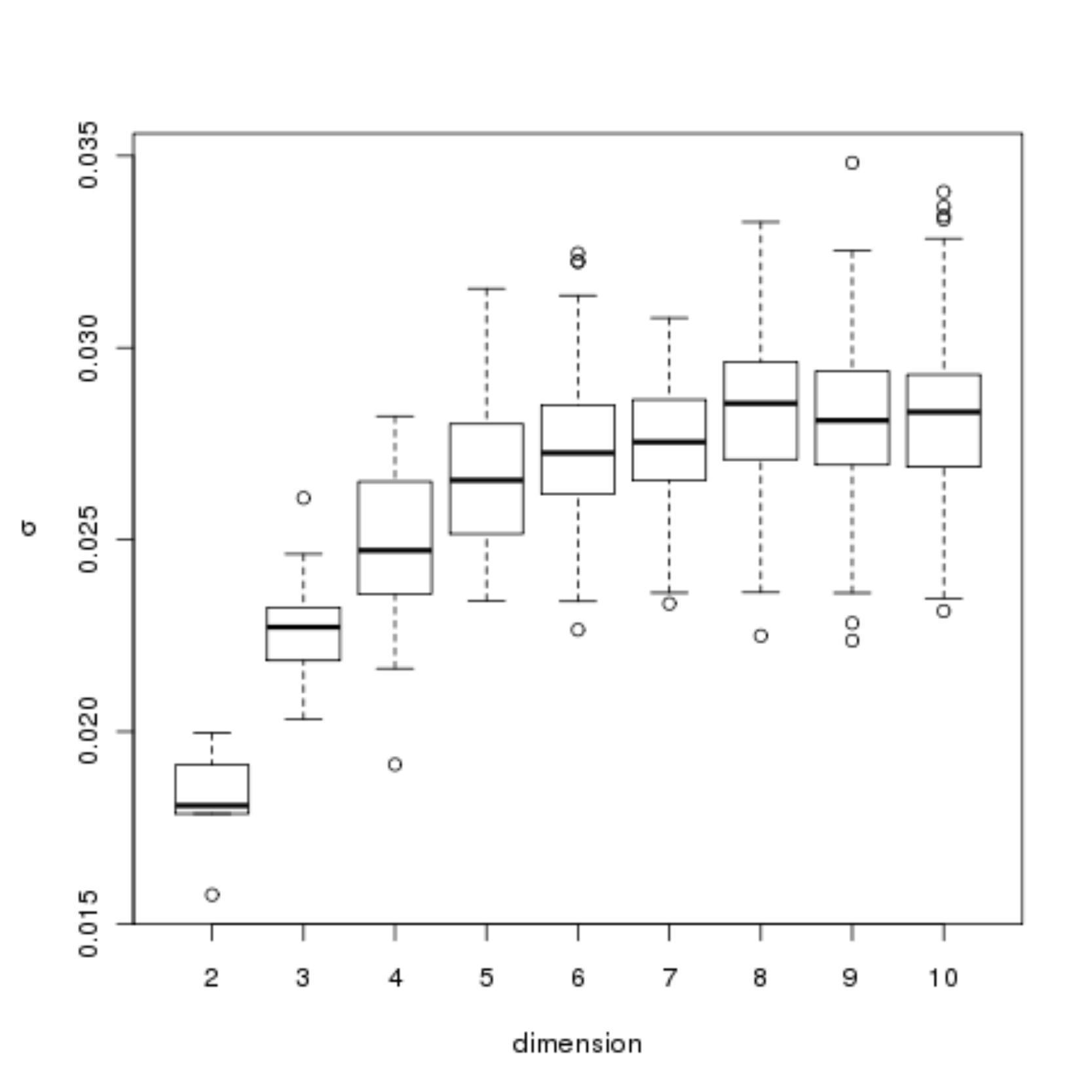}
\end{center}
\vspace*{-1cm}
\caption{MST mean $m$ (left) and standard deviation $\sigma$ (right)
of the 2D subsamples of $C^2$-discrepancy optimized  LHS}
\label{MSTC2}
\end{figure}

Figure~\ref{MSTmaximin} shows
that the MST built over the 2D subsamples of maximin LHS have small
$m$ and high $\sigma$, even for small dimensions $d$,
because of inherent alignments in the entire design. 
Otherwise, in the case LHS optimized with $C^2$-discrepancy, the decrease of $m$
and increase of $\sigma$ for the 2D subsamples are gradual from $d = 1$
to $d = 10$: see Figure~\ref{MSTC2}.
As a consequence, we conclude that $C^2$-discrepancy optimized LHS
of dimension $d \leq 10$ are less subject than maximin designs to the
clustering of points due to projections to 2D subspaces.

\section{Discussion and conclusion}

\subsection{Synthesis of the works}

This paper has considered several issues relative to the practice of the building
of SFD.
Some industrial needs have first been given as challenges: high dimensional
SFD are often required (several tens of variables), while preserving good
space filling properties on the design subprojections.
Several measures of space filling have been discussed:
\Ltwo-discrepancies, mindist interpoint distance criterion (maximin designs),
as well as the recently introduced MST criteria.
It has been shown that the MST criteria are preferable to the
well-known mindist criterion to assess the quality of SFD.
The ESE algorithm is shown by our
numerical tests to behave more efficiently than the MM SA one, since it generally
leads to better improvements of the criterion during the first iterations
(for a same number of elementary perturbations) and eventually provides
similar final values for large computational durations. Thus, the ESE algorithm
appears a relevant choice to carry out LHS optimizations for large dimensions $d$,
when restriction on the computational budget is likely to have an
impact on the result of optimization.
The relevance of the regularized criterion $\phi_p$
to compute maximin LHS has also been numerically checked.

Another contribution of this paper is the deep analysis of the space filling
properties of the two-dimensional subsamples.
Intensive numerical tests have been performed for that purpose, relying on
\Ltwo-discrepancies and minimum spanning tree criteria,
with designs of size $N=100$ and of dimension
$d$ ranging from $2$ to $54$, which is,
to our knowledge, the first numerical study with such a range of design
dimensions.
Among the tested designs (LHS, maximin LHS, several \Ltwo-discrepancy
optimized LHS, Sobol' sequence), only the centered ($C^2$) and wrap-around
($W^2$) discrepancy optimized LHS have shown some strong robustness
to 2D projections in high dimension.
This result numerically confirms the theoretical considerations of
\cite{fanli06}. The other tested designs are not robust
to 2D projections when their dimension is larger than $10$.
Additional tests, not shown here, on other types of designs,
bring the same conclusions.

Table~\ref{plex} summarizes the main properties of LHS which have been
observed.
\begin{table}[!ht]
\centering
\caption{Summary of SFD properties (the second column gives a rank according to the MST
partial order relation).}
\label{plex}
\bigskip
\begin{tabular}{|l|c|c|c|}
\hline
Design & MST & Uniformity on $2$D projections    \\ \hline
LHS &$4$&  no \\ \hline
maximin LHS &$1$ &  no\\  \hline
low C$2$-discrepancy LHS &$2$&  yes \\ \hline
low W$2$-discrepancy LHS &$3$&  yes \\ \hline
\end{tabular}
\end{table}
Eventually, according to our numerical tests, the best LHS are the
$C^2$-discrepancy optimized ones, since 2D subprojections are then quite uniform,
and, in addition, their MST criteria are better than the ones of the $W^2$-discrepancy
optimized LHS.

\subsection{Feedback on the motivating examples}

We illustrate now the benefits
achieved by the use of optimized LHS designs on the two
motivating examples of section~\ref{sec:examples}. 

A typical product of an exploratory study is a metamodel of the
computer code from a first set of inputs-outputs (see introduction).
In \cite{simlin01} and \cite{ioobou10}, several numerical experiments
showed the advantage (smaller mean error) of discrepancy-optimized LHS
rather than SRS or non-optimized LHS to build an acceptable metamodel.
In \cite{marmar12}, a nuclear safety study involving some accidental scenario
simulations was performed using the thermal-hydraulic computer code
of section \ref{sec:nuclear} ($d=31$ input variables).
Because of the high cpu time cost of this code, a metamodel had to be developed.
For this purpose, a $W^2$-discrepancy optimized LHS of size $N=300$ and the
corresponding outputs were computed.
The predictivity coefficient of the fitted Gaussian process metamodel
(see \cite{fanli06}) reached $0.87$, which means that $87\%$ of the model output
variability was explained by the metamodel.
This satisfactory result allowed the authors to use the metamodel
to find an interesting set of inputs so as to estimate low quantiles of the
model output.

In \cite{ioopop12}, a global sensitivity analysis of the
periphyton-grazers sub-model of the prey-predator simulation chain of
section \ref{sec:prey} was carried out. Because of the large number of input
variables, Derivative-based Global Sensitivity Measures (DGSM) were used
\citep{lamioo13}. A DGSM  is defined as the mean of a square model
partial derivative:
$\displaystyle \mbox{DGSM}(X_i) = \ee\left[\left(\frac{\partial G(\bX)}{\partial X_i}\right)^2\right]$.
\cite{kucrod09} showed that low-discrepancy sequences are preferable than SRS to
estimate DGSM.
The partial derivatives of the periphyton-grazers sub-model were
computed by finite differences at each point of a $C^2$-discrepancy optimized LHS of
size $N=100$ and dimension $d=20$.
In total, $N\,(d+1)=2100$ model evaluations were performed.
A study showed that the estimated DGSM had converged, a
result which the SFD directly contributed to.

\subsection{Perspectives}

As a first perspective, analyses of section~\ref{sec:robustness}
could be extended to subprojections of larger dimensions.
It is worth noting that, in a preliminary study, \cite{mar08} confirmed our
conclusions about design subprojection properties by considering 3D
subsamples.

Another future work would be to carry out a more exhaustive and deeper
benchmark of optimization algorithms of LHS.
In particular, an idea would be to look at the convergence of the maximin LHS to
the exact solutions.
These solutions are known in several cases (small $N$ and small $d$) for the
non randomized maximin LHS (see \cite{vanhus07}).

All the SFD of the present work are based on LHS which can be judged as
a rather restrictive design structure.
Building SFD with good properties in the whole space of possible designs is a
grand challenge for future works; see \cite{promul12} for first ideas on this
subject.

Finally, lets us note that a complete emphasis has been
put on the model input space in this paper which deals with
initial exploratory analysis.
Of course, the behaviour of the output variables and the objective of the
study has generally to be considered afterwards with great caution.

\section{Acknowledgments}

Part of this work has been backed by French National Research Agency (ANR)
through the COSINUS program (project COSTA BRAVA noANR-09-COSI-015). 
All our numerical tests were performed within the R statistical software
environment with the DiceDesign package.
We thank Luc Pronzato for helpful discussions and Catalina Ciric for
providing the prey-predator model example.
Finally, we are grateful to both of the reviewers for their valuable comments
and their help with the English.


\bibliographystyle{apalike}
\bibliography{./articles,./books,./thesis,./proceedings,./inbooks,./presentation}

\end{document}